\documentclass[a4paper,12pt]{article}  
\usepackage{graphicx}
\usepackage{amsmath}
\usepackage{amssymb}
\usepackage{enumerate}
\usepackage{amsthm}
\usepackage[T1]{fontenc}
\usepackage[latin1]{inputenc}
\usepackage[round]{natbib}

\newtheorem{thm}{Theorem}
\newtheorem{prop}[thm]{Proposition}
\newtheorem{cor}[thm]{Corollary}
\newtheorem{lemma}[thm]{Lemma}
\theoremstyle{remark}
\newtheorem{remark}[thm]{Remark}
\newtheorem{example}[thm]{Example}
\theoremstyle{definition}

\newcommand{\pr}[2]{P_{#1}\left(#2\right)}

\newcommand{\E}[2]{E_{#1}\left[#2\right]}
\newcommand{\indic}[1]{1_{\left\{#1\right\}}}
\newcommand{\abs}[1]{\left\vert#1\right\vert}

\newcommand{\set}[1]{\left\lbrace #1 \right\rbrace}

\date{\today}

\begin{document}

\title{Survival probabilities of some iterated processes}
\author{\renewcommand{\thefootnote}{\arabic{footnote}}{\sc Christoph Baumgarten}\footnotemark[1]}
\date{\today}

\footnotetext[1]{
Technische Universit\"at Berlin, Institut f\"ur Mathematik, Sekr.\ MA 7-4,
Stra{\ss}e des 17.\ Juni 136, 10623 Berlin, Germany. Email: {\sl baumgart@math.tu-berlin.de}.
}
\maketitle
\begin{abstract}
We study the asymptotic behaviour of the probability that a stochastic process $(Z_t)_{t \geq 0}$ does not exceed a constant barrier up to time $T$ (the so called survival probability) when $Z$ is the composition of two independent processes $(X_t)_{t \in I}$ and $(Y_t)_{t \geq 0}$. To be precise, we consider $(Z_t)_{t \geq 0}$ defined by $Z_t = X \circ \abs{Y_t}$ when $I = [0,\infty)$ and $Z_t = X \circ Y_t$ when $I = \mathbb{R}$.\\
For continuous self-similar processes $(Y_t)_{t \geq 0}$, the rate of decay of survival probability for $Z$ can be inferred directly from the survival probability of $X$ and the index of self-similarity of $Y$. As a corollary, we obtain that the survival probability for iterated Brownian motion decays asymptotically like $T^{-1/2}$. \\
If $Y$ is discontinuous, the range of $Y$ possibly contains gaps which complicates the estimation of the survival probability. We determine the polynomial rate of decay for $X$ being a L\'{e}vy process (possibly two-sided if $I = \mathbb{R}$) and $Y$ being a L\'{e}vy process or random walk under suitable moments conditions.  
\end{abstract}
\newpage
\section{Introduction}
\subsection{Statement of the problem}
The one-sided exit problem consists in finding the asymptotic behaviour of
\begin{equation}\label{eq:one-sided-exit-problem}
   \pr{}{\sup_{t \in [0,T]} Z_t \leq 1}, \quad T \to \infty,
\end{equation}
for a given stochastic processes $Z = (Z_t)_{t \geq 0}$. The probability in \eqref{eq:one-sided-exit-problem} is often called survival or persistence probability up to time $T$. Since it usually cannot be computed explicitly, one aims to specify its asymptotic behaviour. If it decreases polynomially (modulo terms of lower order), i.e.
\[
   \pr{}{\sup_{t \in [0,T]} Z_t \leq 1} = T^{-\theta + o(1)}, \quad T \to \infty,
\]
we call $\theta > 0$ the survival exponent. \\
Of course, \eqref{eq:one-sided-exit-problem} is a classical problem that has been studied for some particular processes such as random walks, Brownian motion with moving boundaries, integrated Brownian motion, fractional Brownian motion (fBm), and other Gaussian processes. Apart from pure theoretical interest, survival probabilities appear in many applications. For instance, the one-sided exit problem arises in various physical models such as reaction diffusion systems and granular media, see the survey of \cite{majumdar:1999} for more examples. Moreover, the study of the one-sided exit problem was motivated by the investigation of the inviscid Burgers equation, see e.g.\ \cite{sinai:1992, bertoin:1998, molchan:1999a}. For a relation to questions about random polynomials and more applications, we refer to \cite{li-shao:2004}. \\
In this article, we consider the one-sided exit problem for processes $Z = (X \circ \abs{Y_t})_{t \geq 0}$ where $X = (X_t)_{t \geq 0}$ and $Y = (Y_t)_{t \geq 0}$ are independent stochastic processes and $Z = (X \circ Y_t)_{t \geq 0}$ if $X = (X_t)_{t \in \mathbb{R}}$ ($\circ$ denotes function composition). 
Such processes will be referred to as \textit{iterated processes}. Starting with the work of \cite{burdzy:1993}, the study of iterated Brownian motion has attracted a lot of interest. Moreover, there are interesting connections of the exit times of iterated processes and the solution of certain fourth-order PDEs (see e.g.\ \cite{allouba-zheng:2001} and \cite{nane:2008}). The asymptotics of the survival probabilities of subordinated Brownian motion is also relevant for the study of Green functions (see e.g.\ \cite{grzywny-ryznar:2008}). However, the one-sided exit problems for itereted processes has not been studied systematically so far. Here we investigate how the survival exponent of $X \circ \abs{Y}$ and $X \circ Y$ is related to that of the outer process $X$ and properties of the inner process $Y$. The relevant scenario affecting the survival probability can be identified so that the results are quite intuitive. For small deviation probabilities (i.e.\ two-sided exit problems), this problem has been investigated by \cite{aurzada-lifshits:2009}. \\
Finally, let us introduce some notation and conventions: If $f,g: \mathbb{R} \to \mathbb{R}$ are two functions, we write $f \precsim g$ $(x \to \infty)$ if $\limsup_{x \to \infty} f(x)/g(x) < \infty$ and $f \asymp g$ if $f \precsim g$ and $g \precsim f$. Moreover, $f \sim g$ $(x \to \infty)$ if $f(x)/g(x) \to 1$ as $x \to \infty$. If $(X_t)_{t \geq 0}$ is a stochastic process, it will often be convenient to write $X(t)$ instead of $X_t$. If $(X_n)_{n \in \mathbb{N}}$ is a discrete time process, we set $X_t = X_{\lfloor t \rfloor}$. Moreover, we say that $(X_t)_{t \in I}$ is self-similar of index $H$ if $(Y_{ct})_{t \geq 0} \stackrel{d}{=} (c^{H} Y_{t})_{t \geq 0}$ for all $c > 0$ where $\stackrel{d}{=}$ denotes equality in distribution.
\subsection{Main results}
First, we consider processes $(X_t)_{t \geq 0}$ and $(Y_t)_{t \geq 0}$ where $Y$ is self-similar and continuous. In this setup, the following result can be established without much difficulty:
\begin{thm}\label{thm:summary_self_sim_1}
   Let $(X_t)_{t \geq 0}$ and $(Y_t)_{t \geq 0}$ be independent stochastic processes. Assume that $Y_0 = 0$ and that $Y$ has continuous paths. Moreover, suppose that $Y$ is self-similar of index $H$. Let $\theta > 0$ and assume that
\[
   \pr{}{\sup_{t \in [0,T]} X_t \leq 1} \asymp T^{-\theta}, \quad T \to \infty, 
\]
and for some $\rho > \theta$,
\begin{equation}\label{eq:small_dev_prop}
   0 < \pr{}{ \sup_{t \in [0,1]} \abs{Y_t} \leq \epsilon } \precsim \epsilon^{\rho}, \quad \epsilon \downarrow 0.
\end{equation}
Then 
\[
 \pr{}{\sup_{t \in [0,T]} X(\abs{Y_t}) \leq 1} \asymp T^{-\theta H}, \quad T \to \infty.
\]
\end{thm}
We remark that the assumption in \eqref{eq:small_dev_prop} (on the so called small deviations of $Y$) is very weak since this probability usually decays exponentially fast. Moreover, the result can be explained quite intuitively: by self-similarity of $Y$, typical fluctuations of $\abs{Y}$ up to time $T$ are of order $T^H$. The rare event that $X$ stays below $1$ until time $T^H$ is then of order $T^{-\theta H}$. The assumption in \eqref{eq:small_dev_prop} prevents a contribution of the event that $Y$ stays close to the origin to the survival exponent of $Z = X \circ \abs{Y}$. In short, the survival probability of $Z$ is determined by a rare event for $X$ and a typical scenario for $Y$. \\
We present various examples in Section~\ref{sec:sup_self_sim_one_sided}. For instance, if $X$ and $Y$ are independent Brownian motions, the survival exponent of $X \circ \abs{Y}$ is $1/4$. \\
The assumption of continuity of the inner process $Y$ allows us to write
\[
   \pr{}{\sup_{t \in [0,T]} X(\abs{Y_t}) \leq 1} = \pr{}{X_t \leq 1, \forall t \in [0,(-I_T) \vee M_T]}
\]
where $I$ and $M$ denote the infimum resp.\ supremum process of $Y$. This will simplify the proof of the upper bound of Theorem~\ref{thm:summary_self_sim_1} very much. If $Y$ is discontinuous, the equality sign has to be replaced by $\geq$ in the preceding equation. It is then a by far more challenging task to find the survival exponent of $X \circ \abs{Y}$ since the gaps in the range of $\abs{Y}$ have to be taken into account. We prove the following theorem for $X$ being a L\'{e}vy process and $Y$ being a random walk or a L\'{e}vy process. 
\begin{thm}\label{thm:summary_LP_at_RW_times}
   Let $(X_t)_{t \geq 0}$ be a L\'{e}vy process such that $\E{}{X_1} = 0$, $\E{}{X_1^2} > 0$ and $\E{}{\exp\left(\abs{X_1}^\alpha\right)} < \infty$ for some $\alpha > 0$. Let $(Y_t)_{t \geq 0}$ denote an independent random walk or L\'{e}vy process with $\E{}{Y_1^2} > 0$ and $\E{}{\exp\left(\abs{Y_1}^\beta \right)} < \infty$ for some $\beta > 0$. 
\begin{enumerate}
   \item If $\E{}{Y_1} = 0$, then $\pr{}{\sup_{t \in [0,T]} X(\abs{Y_t}) \leq 1} = T^{-1/4 + o(1)}.$
\item If $\E{}{Y_1} \neq 0$, then $\pr{}{\sup_{t \in [0,T]} X(\abs{Y_t}) \leq 1} = T^{-1/2 + o(1)}.$
\end{enumerate}
\end{thm}
The lower order terms can be specified more precisely, see Theorem~\ref{thm:LP_at_RW_times} and Theorem~\ref{thm:LP_at_Levy_times}. Again, the results are intuitive: If $\E{}{Y_1} = 0$, the random walk oscillates and typical fluctuations up to time $N$ are of magnitude $\sqrt{N}$. Since the survival exponent $\theta$ of a centered L\'{e}vy process with second finite moments is $1/2$, it is very plausible that the survival exponent of $X \circ \abs{Y}$ is $1/4$ at least if the gaps in the range of the random walk are not too large. If $\E{}{Y_1} > 0$, then $\E{}{Y_N} / N \to \E{}{Y_1}$ by the law of large numbers and one expects the survival exponent of $X \circ \abs{Y}$ to be $1/2$ by the same reasoning. \\
We also exhibit an example showing that an analogous result to Theorem~\ref{thm:summary_LP_at_RW_times} does not hold if the increments of $X$ are not stationary (cf.\ Remark~\ref{rem:counter-ex}) which explains the restriction to L\'{e}vy processes.\\
Up to now, the outer process $X = (X_t)_{t \geq 0}$ had the index set $[0,\infty)$, so it was only possible to evaluate $X$ over the range of the absolute value of the inner process $Y$. In order to consider the one-sided exit problem for $X \circ Y$, we define two-sided processes $X = (X_t)_{t \in \mathbb{R}}$ where
\begin{equation}\label{eq:def-two-sided-proc}
 X_t := 
\begin{cases}
 X^+_t, &\quad t \geq 0,\\
X^-_{-t}, &\quad t <0,
\end{cases}
\end{equation}
and $(X^+)_{t \geq 0}$ and $(X^-_t)_{t \geq 0}$ are independent stochastic processes. We refer to $X^+$ and $X^-$ as the branches of $X$. We prove that the previous results can be extended in a natural way for two-sided processes. 
\begin{thm} \label{thm:summary_self_sim_2}
Let $(X_t)_{t \in \mathbb{R}}$ be a two-sided process generated by $X^-$ and $X^+$ with $X^- \stackrel{d}{=} X^+$. Assume that
\[
	\pr{}{\sup_{t \in [0,T]} X_t \leq 1} \asymp T^{-\theta}
\]
for some $\theta > 0$. Let $(Y_t)_{t \geq 0}$ denote an independent self-similar process of index $H$ with $Y_0 = 0$ and continuous paths such that, as $\epsilon \downarrow 0$,
\begin{align*}
   0 < \pr{}{-\inf_{t \in [0,1]} Y_t \leq \epsilon} \precsim \epsilon^{\eta}, \quad 0 < \pr{}{\sup_{t \in [0,1]} Y_t \leq \epsilon} \precsim \epsilon^{\eta},\quad \pr{}{\sup_{t \in [0,1]} \abs{Y_t} \leq \epsilon} \precsim e^{- \epsilon^{-\gamma}},
\end{align*}
for some $\eta > \theta$ and $\gamma > 0$.
Then
\[
 \pr{}{\sup_{t \in [0,T]} X(Y_t) \leq 1} \asymp T^{-2H\theta}, \quad T \to \infty.
\] 
\end{thm}
We see that the survival exponent in the two-sided setting is twice the exponent of Theorem~\ref{thm:summary_self_sim_1}. This is quite intuitive since by independence of $X^+$ and $X^-$, we have that
\[
   \pr{}{\sup_{t \in [-T,T]} X_t \leq 1} = \pr{}{\sup_{t \in [0,T]} X^+_t \leq 1} \pr{}{\sup_{t \in [0,T]} X^-_t \leq 1} \asymp T^{-2\theta}.
\]
Since the fluctuations of $Y$ up to time $T$ are of magnitude $\pm T^H$ with high probability (this is again ensured by the conditions on $Y$ which are stronger than in Theorem~\ref{thm:summary_self_sim_1}), Theorem~\ref{thm:summary_self_sim_2} appears very natural.\\
In Theorem~\ref{thm:summary_self_sim_2}, we have assumed that the branches of $X$ have the same distribution. This was done for simplicity of exposition, see Theorem~\ref{thm:cont_self_sim_proc_two_sided} for the general case.\\
As a corollary to Theorem~\ref{thm:summary_self_sim_2}, we obtain that the survival exponent of iterated Brownian motion (using the terminology introduced by \cite{burdzy:1993}) is $1/2$.\\
The result corresponding to Theorem~\ref{thm:summary_LP_at_RW_times} in the two-sided setup is
\begin{thm}\label{thm:summary_two_sided_LP_at_RW_times}
   Let $(X_t)_{t \in \mathbb{R}}$ denote a two-sided L\'{e}vy process with branches $X^+, X^-$ such that $\E{}{X_1^\pm} = 0$, $\E{}{(X_1^\pm)^2} > 0$, $\E{}{\exp\left(\abs{X_1^\pm}^\alpha \right)}<\infty$ for some $\alpha > 0$. Let $(Y_t)_{t \geq 0}$ denote another L\'{e}vy process or random walk independent of $X$ with $\E{}{Y_1^2} > 0$ and $\E{}{\exp\left(\abs{Y_1}^\beta \right)} < \infty$ for some $\beta > 0$.  Then
\[
 \pr{}{\sup_{t \in [0,T]} X(Y_t) \leq 1} = T^{-1/2 + o(1)}, \quad T \to \infty.
\]
\end{thm}
Theorem~\ref{thm:summary_two_sided_LP_at_RW_times} shows that the survival exponent is equal to $1/2$ no matter if $\E{}{Y_1} = 0$ or not (in contrast to Theorem~\ref{thm:summary_LP_at_RW_times}, see Remark~\ref{rem:two-sided-LP} for an explanation). \\
We remark that some processes such as fBm are by definition two-sided processes that cannot be written as in \eqref{eq:def-two-sided-proc} since their branches are not independent. We briefly touch upon that case in Section~\ref{sec:fBm}.\\
The remainder of the article is organized as follows. In Section~\ref{sec:sup_self_sim_one_sided}, we assume that the inner process $Y$ is a continuous self-similar process. We compute the survival exponent of $X \circ \abs{Y}$ (Theorem~\ref{thm:summary_self_sim_1}) and provide a couple of examples. Next, we turn to discontinuous processes $Y$. The survial exponent of $X \circ \abs{Y}$ is found for $X$ being a L\'{e}vy process and $Y$ being a random walk or L\'{e}vy process (Theorem~\ref{thm:summary_LP_at_RW_times}) in Section~\ref{sec:discont-one-sided}. Finally, we extend the previous results to two-sided processes (Theorem~\ref{thm:summary_self_sim_2} and Theorem~\ref{thm:summary_two_sided_LP_at_RW_times}) in Section~\ref{sec:two-sided}.
\section{Taking the supremum over the range of a continuous self-similar process}\label{sec:sup_self_sim_one_sided}
If $Y = (Y_t)_{t \geq 0}$ is a stochastic process, denote by $\mathcal{F}^Y_t := \sigma(Y_s : 0 \leq s \leq t)$ the filtration generated by $Y$ up to time $t$. Let us now prove a slightly more general statement than Theorem~\ref{thm:summary_self_sim_1} announced in the introduction. Theorem~\ref{thm:summary_self_sim_1} then follows directly from Theorem~\ref{thm:cont_self_sim_proc} and Lemma~\ref{lem:int_cond_lower_tail_1} below.
\begin{thm}\label{thm:cont_self_sim_proc}
Let $(X_t)_{t \geq 0}$ and $(Y_t)_{t \geq 0}$ be independent stochastic processes. Assume that $Y_0 = 0$ and that $Y$ has continuous paths and is self-similar of index $H$. Let $\theta > 0$.
\begin{enumerate}
 \item If $\pr{}{\sup_{t \in [0,T]} X_t \leq 1} = T^{-\theta + o(1)}$ and if $\E{}{ (\sup_{t \in [0,1]} \abs{Y_t})^{-\eta} } < \infty$ for every $\eta \in (0,\theta)$, then
\[
   \pr{}{\sup_{t \in [0,T]} X(\abs{Y_t}) \leq 1} = T^{-\theta H + o(1)}, \quad T \to \infty.
\]
\item If $\pr{}{\sup_{t \in [0,T]} X_t \leq 1} \asymp T^{-\theta}$ and if $\E{}{ (\sup_{t \in [0,1]} \abs{Y_t})^{-\theta} } < \infty$, then 
\[
 \pr{}{\sup_{t \in [0,T]} X(\abs{Y_t}) \leq 1} \asymp T^{-\theta H}, \quad T \to \infty.
\]
\end{enumerate}
\end{thm}
\begin{proof}
\emph{Upper bound: } Let $\epsilon \in (0,\theta)$. By assumption, we can find constants $C,T_0 > 0$ such that $\pr{}{\sup_{t \in [0,T]} X_t \leq 1} \leq C T^{-\theta + \epsilon}$ for all $T > T_0$. Clearly, we can choose $C$ so large that the inequality holds for all $T > 0$. By continuity of $Y$, the fact that $Y_0 = 0$ and independence of $X$ and $Y$ and self-similarity of $Y$, we have that
   \begin{align*}
       \pr{}{\sup_{t \in [0,T]} X(\abs{Y_t}) \leq 1} &= \E{}{ \pr{}{\sup_{t \in [0,\sup_{u \in [0,T]} \abs{Y_u}]} X_t \leq 1 | \mathcal{F}^Y_T}} \\
	&\leq C \E{}{ \left(\sup_{t \in [0,T]}\abs{Y_t}\right)^{-\theta + \epsilon}} = C \E{}{ \left(\sup_{t \in [0,1]} \abs{Y_t} \right)^{-\theta + \epsilon}} T^{-\theta H + \epsilon H}.
   \end{align*}
Since $(\sup_{t \in [0,1]} \abs{Y_t})^{-\theta + \epsilon}$ is integrable for $\epsilon \in (0,\theta)$, this proves the upper bound in the first case. Under the assumptions of 2., the lines above apply with $\epsilon = 0$.\\
\emph{Lower bound}: Note that 
\begin{align*}
       \pr{}{\sup_{t \in [0,T]} X(\abs{Y_t}) \leq 1} &\geq \pr{}{\sup_{t \in [0,T]} \abs{Y_t} \leq T^H, \sup_{t \in [0,\sup_{u \in [0,T]} \abs{Y_u}]} X_t \leq 1 } \\
&\geq \pr{}{ \sup_{t \in [0,T]} \abs{Y_t} \leq T^H} \, \pr{}{ \sup_{t \in [0,T^H]} X_t \leq 1 } \\
&= \pr{}{ \sup_{t \in [0,1]} Y_t \leq 1} \pr{}{\sup_{t \in [0,T^H]} X_t \leq 1}.
   \end{align*}
This proves the lower bounds.
\end{proof}
\begin{remark}
The proof reveals that the lower bounds of Theorem~\ref{thm:cont_self_sim_proc} are also valid without continuity of paths of $Y$ and the integrability assumption on $Y$. Moreover, the proof of the lower bounds reveals the crucial scenario that determines the survival probability of the composed process. 
\end{remark}
The integrability conditions of Theorem~\ref{thm:cont_self_sim_proc} are satisfied under very mild assumptions on the small deviations of the process $Y$. For convenience, let us state this result in the following lemma.
\begin{lemma}\label{lem:int_cond_lower_tail_1}
Let $Z$ be a random variable such that $Z > 0$ a.s.\ and $\pr{}{Z \leq \epsilon} \precsim \epsilon^\rho$ as $\epsilon \downarrow 0$ for some $\rho > 0$. Then for $\eta \in (0,\rho)$, it holds that
\[
 \E{}{ Z^{-\eta} } < \infty.
\] 
Conversely, if $\E{}{ Z^{-\eta} } < \infty$ for some $\eta >0$, then
\[
  \pr{}{Z \leq \epsilon} \precsim \epsilon^\eta.
\]
\end{lemma}
\begin{proof}
Since $\E{}{Z^{-\eta}} = \E{}{(Z^{-1})^\eta}$ and the latter expectation is finite if $\pr{}{Z^{-1} > x} \precsim x^{-\rho}$ as $x \to \infty$ for some $\rho > \eta$, the first claim follows with $\epsilon = 1/x$.\\
Finally, if $\E{}{ Z^{-\eta} } < \infty$, then for any $\epsilon > 0$, one has
\[
 \E{}{ Z^{-\eta} } \geq \E{}{ Z^{-\eta} ; Z \leq \epsilon } \geq \epsilon^{-\eta} \, \pr{}{Z \leq \epsilon}.
\]
\end{proof}
Although the proof of Theorem~\ref{thm:cont_self_sim_proc} is very simple, the result is applicable to many examples.
\begin{example}\label{ex:BM_iterated}
If $X$ and $Y$ are independent Brownian motions then $\theta = 1/2$ and $H = 1/2$. Since 
\[
 \pr{}{ \sup_{t \in [0,1]} \abs{B_t} \leq \epsilon } \leq C e^{-(\pi^2/8) \, \epsilon^{-2}}, \quad \epsilon > 0,
\]
it is clear that \eqref{eq:small_dev_prop} holds for every $\rho > 0$. Hence, Theorem~\ref{thm:summary_self_sim_1} implies that the survival exponent $X \circ \abs{Y}$ of is $1/4$. \\
More generally, if $W$ and $B^{(1)}, \dots, B^{(n)}$ are independent Brownian motions, it follows for any $n \geq 1$ that 
\[
 \pr{}{ \sup_{t \in [0,T]} W\left(\abs{B^{(1)}} \circ \dots \circ \abs{B^{(n)}_t}\right) \leq 1} \asymp T^{-2^{-(n+1)}}, \quad T \to \infty. 
\] 
\end{example}
\begin{example}
Let $X$ be a process with survival exponent $\theta > 0$. Define
$Y^{(0)}_t = W_t$ where $W$ is a Brownian motion independent of $X$ and define the $n$-times integrated Brownian motion $Y^{(n)}$ for $n \in \mathbb{N}$ recursively by
\[
 Y^{(n)}_t = \int_0^t Y^{(n-1)}_s \, ds, \quad t \geq 0, n \geq 1.
\]
One can check that $Y^{(n)}$ is self-similar with index $H^{(n)} = (2n + 1)/2$. Moreover, the small deviations of $n$-times integrated Brownian motion are known (see Theorem~1.3 of \cite{chen-li:2003}): There exists a constant $\kappa_n \in (0,\infty)$ such that
\[
 \lim_{\epsilon \downarrow 0} \epsilon^{2/(2n+1)} \, \log \pr{}{\sup_{t \in [0,1]} \abs{Y^{(n)}_t} \leq \epsilon } = -\kappa_n, \quad n \geq 1.
\]
In particular, this implies that \eqref{eq:small_dev_prop} is satisfied for any $\rho > 0$. The survival exponent of the iterated process $X \circ \abs{Y^{(n)}}$ is therefore $\theta (2n + 1)/2$ for any $n \geq 1$. In particular, if $X$ is a Brownian motion independent of the Brownian motion $W$, the survival exponent is $(2n + 1)/4$.
\end{example}
\section{Taking the supremum over the range of discontinuous processes}\label{sec:discont-one-sided}
\subsection{Random walks}\label{sec:sup_RW_one_sided}
Let $Y_1, Y_2, \dots$ be a sequence of i.i.d.\ random variables. In the sequel, $S = (S_n)_{n \geq 1}$ denotes the corresponding random walk, i.e.\ $S_n = Y_1 + \dots + Y_n$. Let $\mathcal{F}_N^S := \sigma(S_1,\dots,S_N)$, $M_n := \max_{1 \leq k \leq n} S_k$ and $I_n := \inf_{1 \leq k \leq n} S_k$.\\
The goal of this section is to find the asymptotics of 
\[
  \pr{}{\sup_{n=1,\dots,N} X(\abs{S_n}) \leq 1}, \qquad N \to \infty,
\]
where $X = (X_t)_{t \geq 0}$ is a L\'{e}vy process with $\E{}{X_1} = 0$ and $\E{}{X_1^2} \in (0,\infty)$. First, we recall known results on survival probabilities of L\'{e}vy processes and prove a slight generalization. Under the assumptions on $X$ above, it holds that
\[
   \pr{}{\sup_{t \in [0,T]} X_t \leq 1} \sim c \, T^{-1/2} \, l(T) , \quad T \to \infty,
\]
where $l$ is slowly varying at infinity and $c >0$, see e.g.\ \cite{bingham:1973} or \cite{doney:2007} (Section~4.4) for details. Our goal is to show that the function $l$ may be chosen asymptotically constant which is suggested by the analogous result for random walks: If $(S_n)_{n \geq 1}$ is a centered random walk with finite variance, then $\pr{}{\sup_{n=1,\dots,N} S_n \leq 0} \sim c N^{-1/2}$. However, to the author's knowledge, an analogous result for L\'{e}vy processes has not been stated in the literature so far.\\
Clearly, $\pr{}{\sup_{t \in [0,T]} X_t \leq 1} \leq \pr{}{\sup_{n =1,\dots,\lfloor T \rfloor} X_n \leq 1} \asymp T^{-1/2}$ since $(X_n)_{n \geq 1}$ is a centered random walk with finite variance.  Moreover, if $\E{}{X_1^{2 + \epsilon}} < \infty$ for some $\epsilon > 0$, then also $\pr{}{\sup_{n =1,\dots,\lfloor T \rfloor} X_n \leq 1} \succsim T^{-1/2}$, see e.g.\ Proposition~2.1 in \cite{aurzada-dereich:2011}. The next theorem states the precise asymptotics of $\pr{}{\sup_{t \in [0,T]} X_t \leq 1}$ as $T \to \infty$ under the assumption of finite variance.
\begin{thm}
   Let $(X_t)_{t \geq 0}$ be a L\'{e}vy process such that $\E{}{X_1} = 0$, $\E{}{X_1^2} \in (0,\infty)$. For any $x > 0$, it holds that
\[
   \pr{}{\sup_{t \in [0,T]} X_t \leq x} \sim c(x) \, T^{-1/2}, \quad T \to \infty.
\]

\end{thm}
\begin{proof}
Let $\tau_x$ be the first hitting time of the set $(x,\infty)$, $x > 0$. According to Eq.~4.4.7 of \cite{doney:2007}, it holds that
\begin{equation}\label{eq:LP_tauberian}
   1 - \E{}{e^{-q \tau_x}} \sim U(x) \kappa(q), \quad q \downarrow 0,
\end{equation}
where $U$ is some function (see Eq.~4.4.6 of \cite{doney:2007}) and 
\[
   \kappa(u) = \exp \left( \int_0^\infty \frac{e^{-t} - e^{-ut}}{t} \, \pr{}{X_t > 0} \, dt \right), \quad u \geq 0.
\]
Using that $\int_0^\infty t^{-1}(e^{-t} - e^{-ut}) \, dt = \log u$, it follows that
\[
   \kappa(u) = \sqrt{u} \, \exp \left( \int_0^\infty \frac{e^{-t} - e^{-ut}}{t} \, (\pr{}{X_t > 0} - 1/2) \, dt \right).
\]
We need to show that the integral in the last line converges to a constant as $u \downarrow 0$. To this end, we approximate the term $\pr{}{X_t > 0}$ by $\pr{}{X_n > 0}$ for $t \in [n,n+1]$ which allows us to use classical results from fluctuation theory of random walks to bound the integral from above.\\
Let $u \in (0,1)$ and note that
\begin{align}
   &\int_0^\infty \frac{e^{-ut} - e^{-t}}{t} \, \abs{\pr{}{X_t > 0} - 1/2} \, dt \leq \int_0^1 \frac{1 - e^{-t}}{t} \, \abs{\pr{}{X_t > 0} - 1/2} \, dt \notag \\
&\qquad +\sum_{n=1}^\infty \int_n^{n+1} \frac{e^{-ut} - e^{-t}}{t} \, \abs{\pr{}{X_t > 0} - \pr{}{X_n > 0} } \, dt \notag \\
&\qquad + \sum_{n=1}^\infty \int_n^{n+1} \frac{e^{-ut} - e^{-t}}{t} \, \abs{\pr{}{X_n > 0} - 1/2 } \, dt \notag \\
&\quad \leq c + \sum_{n=1}^\infty n^{-1} \sup_{t \in [n,n+1]} \abs{ \pr{}{X_t > 0} - \pr{}{X_n > 0} } + \sum_{n=1}^\infty n^{-1} \abs{ \pr{}{X_n > 0} - 1/2 }\label{eq:proof_LP_1} 
\end{align}
By a result of \cite{rosen:1962}, it is known that the series $\sum_{n=1}^\infty n^{-1}(\pr{}{X_n > 0} -1/2)$ converges absolutely if $\E{}{X_1} = 0$ and $\E{}{X_1^2} \in (0,\infty)$, so the second series in \eqref{eq:proof_LP_1} converges. Next, we show that the first series also converges. To this end, let $t \in [n,n+1]$. For $f(n) > 0$ (to be chosen appropriately below), we have that
\begin{align*}
   \pr{}{X_t > 0} - \pr{}{X_n > 0} &\leq \pr{}{X_t >0, X_n \leq 0} = \pr{}{X_n \leq 0, (X_t - X_n) + X_n >0} \\
&\leq \pr{}{X_n \leq - f(n), X_t - X_n > f(n)} + \pr{}{0 \geq X_n > -f(n)} \\
&\leq \sup_{u \in [0,1]} \pr{}{X_u > f(n)} + \pr{}{\abs{X_n} \leq f(n)} \\
&\leq \frac{\E{}{X_1^2}}{f(n)^2} + \pr{}{\abs{X_n} \leq f(n)}.
\end{align*}
We have used the independence and stationarity of increments of $X$ and the fact that $\E{}{X_t^2} = t \cdot \E{}{X_1^2}$ for $t \geq 0$ in the above estimates. By the same argument, one shows that $\pr{}{X_n > 0} - \pr{}{X_t > 0} \leq \frac{\E{}{X_1^2}}{f(n)^2} + \pr{}{\abs{X_n} \leq f(n)}$, so for $f(n) = n^{1/6}$, we obtain that
\begin{align*}
   \sum_{n=1}^\infty n^{-1} \sup_{t \in [n,n+1]} \abs{ \pr{}{X_t > 0} - \pr{}{X_n > 0} } \leq \E{}{X_1^2} \sum_{n=1}^\infty n^{-4/3} + \sum_{n=1}^\infty  n^{-1} \pr{}{\abs{X_n} \leq n^{1/6}}.
\end{align*}
We claim that 
\begin{equation}\label{eq:proof_LP_2}
  \pr{}{\abs{X_n} \leq n^{1/6}} = \pr{}{\abs{X_n}/\sqrt{n} \leq n^{-1/3}} \asymp n^{-1/3}, \quad n \to \infty
\end{equation}
If this holds, the first series in \eqref{eq:proof_LP_1} is finite and we conclude that 
\[
    \int_0^\infty \frac{1 - e^{-t}}{t} \, \abs{\pr{}{X_t > 0} - 1/2} \, dt = \sup_{u \in (0,1)} \int_0^\infty \frac{e^{-ut} - e^{-t}}{t} \, \abs{\pr{}{X_t > 0} - 1/2} \, dt < \infty.
\]
Hence,
\[
   \kappa(u) \sim \sqrt{u} \, \exp \left( \int_0^\infty \frac{1 - e^{-t}}{t} \, (\pr{}{X_t > 0} - 1/2) \, dt \right), \quad u \downarrow 0.
\]
The theorem now follows from \eqref{eq:LP_tauberian} by standard Tauberian arguments.\\
It remains to show that \eqref{eq:proof_LP_2} holds which follows from local limit theorems. Indeed, if $X_1$ has a nonlattice distribution (i.e.\ $\abs{ \E{}{\exp(iuX_1)} } < 1$ for $u \in \mathbb{R}\setminus\set{0}$), Theorem~1 of \cite{stone:1965} yields that
\begin{align*}
   \pr{}{X_n/\sqrt{n} \in (x,x+h]} = \pr{}{Z \in (x,x+h]} + o_n(1)(h + n^{-1/2}), \quad n \to \infty,
\end{align*}
where $Z$ is standard normal and the $o(\cdot)$-term is uniform in $h$ and $x$. Hence,
\[
   \pr{}{X_n/\sqrt{n} \in (-n^{-1/3},n^{-1/3}]} \precsim n^{-1/3}, \quad N \to \infty,
\]
since $\pr{}{Z \in [-\epsilon,\epsilon]} \sim \sqrt{2/\pi} \, \epsilon$ as $\epsilon \downarrow 0$. \\
If $X_1$ has a lattice distribution, then \eqref{eq:proof_LP_2} follows from a theorem of \cite{gnedenko-kolmogorov:1968} (p.233).
\end{proof}
Having determined the asymptotics of the survival probability for $X$, let us continue to give some heuristics concerning the survival exponent of $X \circ \abs{S}$. If $\E{}{Y_1} = 0$ and $\E{}{Y_1^2} = 1$, it follows from the invariance principle that
\[
  \lim_{N \to \infty} \pr{}{\sup_{n=1,\dots,N} \abs{S_n} \leq \sqrt{N} \, x } = \pr{}{\sup_{t \in [0,1]} \abs{B_t} \leq x}, \quad x > 0.
\]
Here, $B$ denotes a standard Brownian motion. Intuitively, one would therefore expect that
\[
    \pr{}{\sup_{n=1,\dots,N} X(\abs{S_n}) \leq 1} \asymp  \pr{}{\sup_{t \in [0,\sqrt{N}]} X_t \leq 1} \asymp N^{-1/4},
\]
at least if the points $\abs{S_1},\dots,\abs{S_N}$ are sufficiently ``dense`` in $[0,\sqrt{N}]$. Under suitable moment conditions on the random walk, we show that the survival exponent is indeed $1/4$. For simplicity of notation, we denote by $\mathcal{X}(\gamma)$ the class of random variables $X$ with $\E{}{X^2} > 0$ and $\E{}{e^{\abs{X}^\gamma}} < \infty$ where $\gamma > 0$.
\begin{thm}\label{thm:LP_at_RW_times}
   Let $(X_t)_{t \geq 0}$ denote a L\'{e}vy process with $\E{}{X_1} = 0$ and $X_1 \in \mathcal{X}(\alpha)$ for some $\alpha \in (0,1]$. Assume that $(Y_n)_{n \geq 1}$ is independent of $X$ with $Y_1 \in \mathcal{X}(\beta)$ for some $\beta \in (0,1]$. 
\begin{enumerate}
 \item If $\E{}{Y_1} = 0$, then
  \[
    N^{-1/4} \precsim \pr{}{\sup_{n=1,\dots,N} X(\abs{S_n}) \leq 1} \precsim N^{-1/4} \, (\log N)^{1/4 + 1/(\alpha \wedge \beta)}, \quad N \to \infty.
  \]
  \item 
If $\E{}{Y_1} \neq 0$, then
  \[
    N^{-1/2} \precsim \pr{}{\sup_{n=1,\dots,N} X(\abs{S_n}) \leq 1} \precsim N^{-1/2} \, (\log N)^{1/(\alpha \wedge \beta)}, \quad N \to \infty.
  \]
\end{enumerate}
In either case, the lower bound also holds without the assumption of stretched exponential moments on $X_1$ and $Y_1$.
\end{thm}
The lower bounds follow from the lowers bounds of Theorem~\ref{thm:LP_at_Levy_times} below in view of \eqref{eq:comp_LP_vs_RW}.\\
\emph{Upper bound:}
Let us first introduce some more notation: Denote by $\sigma(n)$ the $n$-th time that the random walk $S$ reaches a new maximum. Then $M_{\sigma(n)} = S_{\sigma(n)}$ is the position of the random walk at that time and
\[
   S_{\sigma(n)} = \mathcal{H}_1 + \dots + \mathcal{H}_n
\]
where $\mathcal{H}_n = S_{\sigma(n)} - S_{\sigma(n-1)}$ is the $n$-th ladder height and $\mathcal{H}_1, \mathcal{H}_2, \dots$ are i.i.d.\ (see e.g.\ \cite{feller-vol2-1970}, Chapter XII). \\
Before proving the upper bound of Theorem~\ref{thm:LP_at_RW_times}, we need two auxiliary results.
\begin{lemma}\label{lem:sup_LP_ladder_height}
Under the assumptions of Theorem~\ref{thm:LP_at_RW_times}, if $\E{}{Y_1} \geq 0$, then for any $p >0$, one can find a constant $c>0$ such that
\[
 \pr{}{\sup_{t \in [0,\mathcal{H}_1]} X_t \geq c (\log N)^{1/(\alpha \wedge \beta)} } \precsim N^{-p}, \quad N \to \infty.
\]
\end{lemma}
\begin{proof}
Note that
\begin{align}
   & \pr{}{ \sup_{t \in [0,\mathcal{H}_1]} X_t > c (\log N)^{1/(\alpha \wedge \beta)} } \notag \\
&\quad \leq \pr{}{\sup_{t \in [0,\mathcal{H}_1]} X_t > c (\log N)^{1/(\alpha \wedge \beta)}, \mathcal{H}_1 \leq d (\log N)^{1/\beta}} + \pr{}{\mathcal{H}_1 > d (\log N)^{1/\beta}} \notag \\
& \quad \leq  \pr{}{\sup_{t \in [0,d (\log N)^{1/\beta}]} X_t > c (\log N)^{1/(\alpha \wedge \beta)}} + \pr{}{\mathcal{H}_1 > d (\log N)^{1/\beta}}. \label{eq:proof_ladder_height}
\end{align}
We distinguish the cases $\E{}{Y_1} = 0$ and $\E{}{Y_1} > 0$. \\
\emph{Case 1}: First, assume that $\E{}{Y_1} = 0$. The second term in \eqref{eq:proof_ladder_height} can be controlled by Chebychev's inequality and a result on the moments of ladder heights of \cite{doney:1980}: Set $\varphi(x) =  e^{x^\beta} x^{\beta - 1}$ on $[x_0,\infty)$ where $x_0$ is chosen in such a way that $\varphi$ is increasing on this interval. On $[0,x_0]$, let $\varphi$ be a non-negative bounded increasing function such that $\varphi$ is differentiable and increasing on $(0,\infty)$. Set $\Phi(x) := \int_0^x \varphi(u) \, du$, $x \geq 0$. Then $\Phi$ is bounded on $[0,x_0]$ by some constant $C > 0$ and for $x > x_0$, we have 
\[
 \Phi(x) \leq C + \int_{x_0}^x e^{u^\beta} u^{\beta - 1} \, du = C + \beta^{-1} \left( e^{x^\beta} - e^{x_0^\beta}\right) \leq C + \beta^{-1} \, e^{x^\beta}.
\]
Therefore,
\[
  \E{}{\Phi(\abs{Y_1})} \leq 2C + \beta^{-1} \, \E{}{e^{\abs{X_1}^\beta} \indic{\abs{Y_1} \geq x_0}} < \infty
\]
since $Y_1 \in \mathcal{X}(\beta)$. By Theorem~1 of \cite{doney:1980}, this implies $\E{}{\varphi(\mathcal{H}_1)} < \infty$. In particular, for large $N$,
\begin{align*}
  \pr{}{\mathcal{H}_1 > d (\log N)^{1/\beta}} \leq \frac{\E{}{\varphi(\mathcal{H}_1) }}{\varphi(d (\log N)^{1/\beta})} = \frac{\E{}{\varphi(\mathcal{H}_1)}}{ d^{\beta - 1}(\log N)^{(\beta - 1)/\beta}} \, \exp(- d^\beta \, \log N).
\end{align*}
It is clear that this term is $o(N^{-p})$ if $d$ is sufficiently large.\\
Let us now consider the first term in \eqref{eq:proof_ladder_height}. Assume first that $\alpha = 1$. Since $X$ is a L\'{e}vy process, $\E{}{e^{X_t}} = e^{\Lambda t}$ for all $t \geq 0$ and $\Lambda = \log \E{}{e^{X_1}}$. Moreover, $(\exp(X_t))_{t \geq 0}$ is a positive submartingale since $X$ is a martingale and therefore, it follows from Doob's inequality that
\[
 \pr{}{\sup_{t \in [0,T]} X_t > x} \leq e^{- x} \E{}{e^{X_T}} = e^{-x + \Lambda T}, \quad T,x >0.
\]
We therefore get for the first term in \eqref{eq:proof_ladder_height} that
\begin{align*}
   \pr{}{\sup_{t \in [0,d (\log N)^{1/\beta}]} X_t > c (\log N)^{1/(\alpha \wedge \beta)},} \leq e^{-c (\log N)^{1/\alpha \wedge \beta} + \Lambda d (\log T)^{1/\beta} } \precsim N^{-p}
\end{align*}
if $c$ is chosen large enough since $\alpha \wedge \beta = \beta \leq 1$. \\
If $\alpha < 1$, we apply Cr\'{a}mer's theorem without exponential moments. To this end, recall the following maximal inequality for L\'{e}vy processes:
\begin{equation}\label{eq:maximal_ineq_LP}
   \pr{}{\sup_{t \in [0,T]} \abs{X_t} \geq x} \leq 9 \, \pr{}{X_T \geq x/30}, \quad T,x > 0.
\end{equation}
This follows from Montgomery-Smith's inequality for sums of centered i.i.d.\ random variables (Corollary 4 of \cite{montg-smith:1993}) since 
\begin{align*}
   \pr{}{\sup_{t \in [0,T]} \abs{X_t} \geq x} = \lim_{n \to \infty} \pr{}{\sup_{k=1,\dots,n} \abs{X_{k T/n}} \geq x} \leq 9 \pr{}{\abs{X_T} \geq x/30}.
\end{align*}
The application of Montgomery-Smith's inequality is possible since $X_{k T/n} = Y_{1,T/n} + \dots + Y_{k,T/n}$ ($k = 1,\dots,n)$ where $Y_{1,T/n},\dots,Y_{n,T/n}$ are i.i.d.\ random variables with $Y_{1,T/n} \stackrel{d}{=} X_{T/n}$.\\
Let $S = (S_n)_{ n \geq 1}$ denote a random walk whose increments have the same law as $X_1$. Then we deduce from \eqref{eq:maximal_ineq_LP} that
\begin{equation}\label{eq:proof_ladder_height_2}
   \pr{}{\sup_{t \in [0,T]} X_t \geq x} \leq 9 \pr{}{\abs{S_{\lceil T \rceil}} \geq x/30} = 9 \left( \pr{}{S_{\lceil T \rceil} \geq x/30} + \pr{}{-S_{\lceil T \rceil} \geq x/30} \right).
\end{equation}
Next, it suffices to apply a large deviations result under the assumption of stretched exponential moments (see Eq.\ 2.32 in \cite{nagaev:1979}): There is a constant $C_1 > 0$ such that for $N$ and $x > 0$, one has that
\begin{align*}
   \pr{}{S_N > \sigma x} &\leq C_1 \left(e^{-\sigma^2 x^2/(20 N)} + N \, \pr{}{X_1 > \sigma x/2}\right) \\
&\leq C_1 \left(e^{-\sigma^2 x^2/(20 N)} + N e^{-(\sigma x/2)^\alpha} \, \E{}{e^{\abs{X_1}^\alpha}}\right),
\end{align*}
where $\sigma^2 := \E{}{X_1^2}$. Hence, combining this with \eqref{eq:proof_ladder_height_2}, we have for some $c$ and all $N$ large enough that
\[
   \pr{}{\sup_{t \in [0,d (\log N)^{1/\beta}]} X_t \geq c (\log N)^{1/(\alpha \wedge \beta)} } \precsim N^{-p}.
\]
\emph{Case 2:} Assume that $\E{}{Y_1} > 0$. The tail behaviour of the first ladder height $\mathcal{H}_1$ can be determined in view of the following estimates: For $x > 0$, one has
\begin{align*}
 \pr{}{Y_1 \geq x} \leq \pr{}{\mathcal{H}_1 \geq x } &= \pr{}{Y_1 \geq x} + \sum_{n=2}^\infty \pr{}{S_1 \leq 0, \dots, S_{n-1} \leq 0, S_{n-1} + Y_n \geq x} \\
&\leq \pr{}{Y_1 \geq x} + \sum_{n=2}^\infty \pr{}{S_1 \leq 0, \dots, S_{n-1} \leq 0, Y_n \geq x} \\
&= \pr{}{Y_1 \geq x} \sum_{n=1}^\infty \pr{}{M_{n-1} \leq 0}.
\end{align*}
It is not hard to check that the latter series converges to a finite value $d_1$ since $\E{}{Y_1} > 0$. Therefore, for any $p >0$, we have that
\[
 \pr{}{\mathcal{H}_1 \geq d(\log N)^{1/\beta} } \leq d_1 \pr{}{Y_1 \geq d (\log N)^{1/\beta} } \leq d_1 \, \frac{ \E{}{e^{\abs{Y_1}^\beta}} }{ \exp( d^\beta \, \log N)} = o(N^{-p}).
\]
for $d$ sufficiently large. Hence, the same arguments used in the first case complete the proof.
\end{proof}

Here and later, we also need the following auxiliary result similar to Proposition~2.1 of \cite{aurzada-dereich:2011}.
\begin{lemma}\label{lem:one-sided-exit-normalized-RW}
Let  $f \colon [0,\infty) \to (0,\infty)$ be a measurable function such that $f(N)/\sqrt{N} \to 0$ as $N \to \infty$. Let $(Y_n)_{n \geq 1}$ denote a sequence of i.i.d.\ random variables with $\E{}{Y_1} = 0$ and $Y_1 \in \mathcal{X}(\beta)$ for some $\beta \in (0,1]$. Let $(S_n)_{ n \geq 1}$ denote the corresponding random walk. 
\begin{enumerate}
 \item If $f(N) \precsim (\log N)^{1/\beta}$, then
\[
 \pr{}{ M_N \leq f(N) } \precsim (\log N)^{1/\beta} \, N^{-1/2} , \quad N \to \infty.
\]
\item If $(\log N)^{1/\beta}/f(N) \to 0$ as $N \to \infty$, then
\[
  \pr{}{ M_N \leq f(N) } \sim \sqrt{\frac{2}{\pi \E{}{Y_1^2} }} \, \frac{f(N)}{\sqrt{N}}, \quad N \to \infty.
\]
\end{enumerate}
\end{lemma}
\begin{proof} 
Let $\sigma^2 := \E{}{Y_1^2}$. We need a result on the speed of convergence in the invariance principle under the assumption of stretched exponential moments. According to \cite{sawyer:1968} (p.\ 363, Eq.\ 1.5), it holds that
\begin{equation}\label{eq:sawyer-estimate}
   \sup_{x \geq 0} \abs{\pr{}{\sup_{n=1,\dots,N} S_n \leq x \, \sqrt{\sigma^2 N}} - \pr{}{\abs{B_1} \leq x} } \precsim (\log N)^{1/\beta} \, N^{-1/2}, \quad N \to \infty.
\end{equation}
Hence, since $\sup_{t \in [0,1]} B_t \stackrel{d}{=} \abs{B_1}$, we conclude that
\begin{align}
  \pr{}{ M_N \leq f(N) } &\leq \abs{ \pr{}{M_N  \leq f(N) } - \pr{}{\sup_{t \in [0,1]} B_t \leq f(N)/\sqrt{\sigma^2 N}} } \notag \\
&\quad + \pr{}{\sup_{t \in [0,1]} B_t \leq f(N)/\sqrt{\sigma^2 N} } \notag \\
&\leq C (\log N)^{1/\beta} N^{-1/2} + \sqrt{ 2/(\pi \sigma^2) } \, f(N) \,N^{-1/2}. \label{eq:sawyer-estimate-2}
\end{align}
Depending on the behaviour of $f(N)$ stated in the lemma, the order of the first or second term is the dominant one. \\
The same argument applies for the proof of the lower bound in the second case.
\end{proof}
\begin{remark}\label{rem:sawyer}
Note that due to the uniform estimate in \eqref{eq:sawyer-estimate}, the constant $C$ in \eqref{eq:sawyer-estimate-2} only depends on $N$, but not on the function $f$. This observation will be relevant later on.
\end{remark}
\begin{remark}\label{rem:LP_below_log}
We frequently need to apply Lemma~\ref{lem:one-sided-exit-normalized-RW} to L\'{e}vy processes in the following situation: Let $(X_t)_{t \geq 0}$ denote a L\'{e}vy process such that $\E{}{X_1} = 0$ and $X_1 \in \mathcal{X}(\alpha)$ and let $g \colon [0,\infty) \to (0,\infty)$ be a function such that $g(T) \to \infty$ as $T \to \infty$. For any $c, \rho > 0$, it follows for $T$ large enough that
\begin{align*}
   \pr{}{\sup_{t \in [0,g(T)]} X_t \leq c (\log T)^\rho} &\leq \pr{}{\sup_{n =1,\dots,\lfloor g(T) \rfloor} X_n \leq c (\log T)^\rho} \\
&\leq C \, \frac{(\log g(T) )^{1/\alpha} +  c(\log T)^\rho}{\sqrt{g(T)}}.
\end{align*}
This follows directly from the proof of Lemma~\ref{lem:one-sided-exit-normalized-RW} since $(X_n)_{ n \geq 1}$ is a random walk.
\end{remark}
We are now ready to establish the upper bounds of Theorem~\ref{thm:LP_at_RW_times}.
\begin{proof}(Upper bound of Theorem~\ref{thm:LP_at_RW_times} if $\E{}{Y_1} = 0$.)\\
For the upper bound, the idea is to consider the supremum of the L\'{e}vy process $X$ only at those points where the random walks either reaches a new maximum or a new minimum. More specifically, we have that
   \begin{align*}
      &\pr{}{\sup_{n=1,\dots,N} X(\abs{S_n}) \leq 1} \\
&\quad \leq \pr{}{\max_{n =1,\dots,N} \abs{S_n} \leq \sqrt{N} / f(N)} + \pr{}{M_N \geq \sqrt{N} / f(N), \, \sup_{n=1,\dots,N} X(M_n) \leq 1 } \\
&\qquad + \pr{}{I_N \leq - \sqrt{N} / f(N), \, \sup_{n=1,\dots,N} X(-I_n) \leq 1 } =: J_1(N) + J_2(N) + J_3(N),
   \end{align*}
where $f(N) := \sqrt{\log N }$.\\
\emph{Estimate for $J_1$:}\\
This term is of lower order than $N^{-1/4}$ by a small deviations results of \cite{acosta:1983}. Indeed, by Theorem~4.3 of \cite{acosta:1983}, one has
\begin{equation}\label{eq:de-acosta_result}
 \limsup_{N \to \infty} \frac{\log \pr{}{\sup_{n=1,\dots,N} \abs{S_n} \leq \sqrt{N} / a_N} }{a_N^2} \leq - \pi^2/8
\end{equation}
whenever $0 < a_N \to \infty$ and $a_N^2 / N \to 0$.
This shows that
\begin{equation}\label{eq:proof_RW_term1}
 \pr{}{\max_{n=1,\dots,N} \abs{S_n} \leq \sqrt{N}  / f(N) } = o(N^{-1/4}).
\end{equation}
\emph{Estimate for $J_2$:}\\
Let $A_N := \set{M_N \geq \sqrt{N} / f(N)}$. Conditioning on $\mathcal{F}_N^S$, we get
\[
 J_2(N) = \pr{}{A_N, \, \sup_{n =1,\dots, N} X(M_n) \leq 1 } = \E{}{\indic{A_N} \, \pr{}{ \sup_{n \, : \, \sigma(n) \leq N} X(S_{\sigma(n)}) \leq 1 \vert \mathcal{F}_N^S } }.
\]
We now estimate the term under the expectation sign. For $\rho > 0$ to be specified later, we have
\begin{align*}
   &\indic{A_N} \pr{}{\sup_{n : \sigma(n) \leq N} X(S_{\sigma(n)}) \leq 1 \vert \mathcal{F}_N^S} \leq \indic{A_N} \pr{}{ \sup_{t \in [0,M_N]} X_t \leq 1 + c (\log N)^\rho  \vert \mathcal{F}_N^S} + \\
&\qquad \indic{A_N}  \pr{}{ \bigcup_{n \, :\, \sigma(n) \leq N}  \set{\sup_{t \in [S_{\sigma(n-1)},S_{\sigma(n)}]} X_t - X_{S_{\sigma(n-1)}} > c (\log N)^\rho}  \vert \mathcal{F}_N^S} \\
&\quad \leq \pr{}{ \sup_{t \in [0,\sqrt{N}/f(N)]} X_t \leq 1 + c (\log N)^\rho} \\
&\qquad + \sum_{n=1}^N \pr{}{ \sup_{t \in [0,S_{\sigma(n)} - S_{\sigma(n-1)}]} X_{t + S_{\sigma(n-1)}} - X_{S_{\sigma(n-1)}} > c (\log N)^\rho  \vert \mathcal{F}_N^S}.
\end{align*}
Since $X_1 \in \mathcal{X}(\alpha)$, it follows from Remark~\ref{rem:LP_below_log} that for $N$ sufficiently large and $\rho \geq 1/\alpha$, one has 
\begin{equation}
  \pr{}{ \sup_{t \in [0,\sqrt{N}/f(N)]} X_t \leq 1 + c (\log N)^\rho} \leq  c_1 (\log N)^{\rho + 1/4} \, N^{-1/4}.\label{eq:proof_RW_term2-1}
\end{equation}
Next, since $X$ is a L\'{e}vy process independent of $S$, we have that
\begin{align*}
   &\pr{}{ \sup_{t \in [0,S_{\sigma(n)} - S_{\sigma(n-1)}]} X_{t + S_{\sigma(n-1)}} - X_{S_{\sigma(n-1)}} > c (\log N)^\rho  \vert \mathcal{F}_N^S} \\
&\quad = \pr{}{\sup_{t \in [0,\mathcal{H}_n]} X_t > c (\log N)^\rho  \vert \mathcal{F}_N^S}.
\end{align*}
Using \eqref{eq:proof_RW_term2-1} and keeping in mind that the $\mathcal{H}_n$ are i.i.d., the above estimates imply that
\begin{equation}
  J_2(N) \leq c_1 (\log N)^{\rho + 1/4} \, N^{-1/4} + N \, \pr{}{ \sup_{t \in [0,\mathcal{H}_1]} X_t > c (\log N)^\rho }.
\end{equation}
In view of Lemma~\ref{lem:sup_LP_ladder_height}, for $\rho := 1/(\alpha \wedge \beta)$ and $c$ large enough, we conclude that
\begin{equation}\label{eq:proof_RW_term2}
 J_2(N) \leq c_2 (\log N)^{1/4 + 1/(\alpha \wedge \beta)} \, N^{-1/4} +  o\left(N^{-1/4}\right).
\end{equation}
\emph{Estimate for $J_3$:} \\
Using this time descending ladder epochs and heights (or considering the random walk $(-S_n)_{n \geq 1}$ in the previous step), one can prove analogously that
 \begin{equation}\label{eq:proof_RW_term3}
 J_3(N) \leq c_2 (\log N)^{1/4 + 1/(\alpha \wedge \beta)} \, N^{-1/4} +  o\left(N^{-1/4}\right).
\end{equation}
Combining \eqref{eq:proof_RW_term1}, \eqref{eq:proof_RW_term2} and \eqref{eq:proof_RW_term3} finishes the proof of the upper bound if $\E{}{Y_1} = 0$.
\end{proof}
\begin{proof}(Upper bound of Theorem~\ref{thm:LP_at_RW_times} if $\E{}{Y_1} \neq 0$.)\\
It suffices to prove the lemma for the case $m := \E{}{Y_1} > 0$. The result for $\E{}{Y_1} < 0$  then follows by considering $-Y_1,-Y_2,\dots$.\\
Clearly, we can write
\begin{align}\label{eq:upper_bound_RW_pos_exp}
\pr{}{\sup_{n = 1,\dots,N} X(\abs{S_n}) \leq 1} &\leq \pr{}{M_N \leq c N} + \pr{}{M_N \geq c N, \, \sup_{n=1,\dots,N} X(M_n) \leq 1} \\
\notag &=: J_1(N) + J_2(N). 
\end{align}
\emph{Estimate for $J_1$:}\\
Note that $M_N = \max_{n=1,\dots,N} \left(\tilde{S}_n + m \, n \right) \geq \tilde{S}_N + m \, N$ where $\tilde{S}_n = (Y_1 - m) + \dots + (Y_n - m)$ is a centered random walk. Hence, for $c < m$, one has 
\[
  J_1(N) = \pr{}{M_N \leq c N} \leq \pr{}{\tilde{S}_N \leq N(c - m)} = o(N^{-1/2}), \quad N \to \infty.
\]
\emph{Estimate for $J_2$:}\\
Let $A_N := \set{M_N \geq c N}$. Denote again by $\sigma(n)$ the $n$-th time that the random walk $S$ reaches a new maximum and by $\mathcal{H}_n = S_{\sigma(n)} - S_{\sigma(n-1)}$ the $n$-th ladder height of $S$. Using that the $\mathcal{H}_n$ are i.i.d., $J_2$ can be estimated as above:
\begin{align}
J_2(N) &= \E{}{\indic{A_N} \pr{}{ \sup_{n=1,\dots,N} X(M_n) \leq 1 \vert \mathcal{F}_N^S} } \notag \\
&\leq \E{}{\indic{A_N} \pr{}{ \sup_{t \in [0,M_N]} X_t \leq 1 + c_1 (\log N)^\rho \vert \mathcal{F}_N^S }} \notag \\
& \quad + \E{}{\indic{A_N} \pr{}{ \bigcup_{n \, :\, \sigma(n) \leq N}  \set{\sup_{t \in [S_{\sigma(n-1)},S_{\sigma(n)}]} X_t - X_{S_{\sigma(n-1)}} > c_1 (\log N)^\rho}  \vert \mathcal{F}_N^S} } \notag \\
&\leq \pr{}{ \sup_{t \in [0,c N]} X_t \leq 1 + c_1 (\log N)^\rho } \notag \\
& \quad + \E{}{\indic{A_N} \sum_{n=1}^N \pr{}{ \sup_{t \in [0,\mathcal{H}_n]} X_{t + S_{\sigma(n-1)}} - X_{S_{\sigma(n-1)}} > c_1 (\log N)^\rho  \vert \mathcal{F}_N^S} } \notag \\ 
&\leq c_2 \frac{(\log N)^\rho }{ \sqrt{c N} } + N \pr{}{ \sup_{t \in [0,\mathcal{H}_1]} X_t > c_1 (\log N)^\rho }. \label{eq:sup_discr_to_cont}
\end{align}
The last inequality holds for $N$ sufficiently large and $\rho \geq 1/\alpha$ by Remark~\ref{rem:LP_below_log}. Applying Lemma~\ref{lem:sup_LP_ladder_height}, we conclude that for $\rho := 1/(\alpha \wedge \beta)$ and $c_1$ large enough, we have that
\[
  J_2(N) \precsim (\log N)^{1/(\alpha \wedge \beta)} \, N^{-1/2}, \quad N \to \infty.
\]
\end{proof}
\begin{remark}\label{rem:counter-ex}
One might wonder if the assumption that the outer process $X$ is a L\'{e}vy process can be relaxed. In view of Theorem~\ref{thm:cont_self_sim_proc}, one might guess that if $X$ has a survival exponent $\theta > 0$, it would follow that
\[
   \pr{}{\sup_{n=1,\dots,N} X(\abs{S_n}) \leq 1} = N^{-\theta/2 + o(1)}
\]
under suitable moment conditions. However, this turns out to be false in general. As an example, consider a sequence $\tilde{X}_1, \tilde{X}_2, \dots$ of independent random variables with $\pr{}{\tilde{X}_n = 2} = 1 - \pr{}{\tilde{X}_n = 0} = 1/(n+1)$ for $n \geq 1$ and define $X = (X_t)_{t \geq 0}$ by
\[
   X_t = \tilde{X}_n \quad \text{if } t = (2n - 1)/2 \text{ for some } n \in \mathbb{N}, \quad X_t = 0 \quad \text{else}.
\]
Obviously, $X$ does not have stationary increments. Moreover, it is not hard to check that
\[
   \pr{}{\sup_{t \in [0,T]} X_t \leq 1} \asymp \pr{}{\tilde{X}_1 = 0, \dots, \tilde{X}_{\lfloor T \rfloor} = 0} = \prod_{n=1}^{\lfloor T \rfloor} (1 - 1/(n+1)) \asymp T^{-1}.
\]
If $(S_n)_{n \geq 1}$ is a symmetric simple random walk, one has by construction that $X(\abs{S_n}) = 0$ for all $n$. \\
If $X$ has stationary, but not necessarily independent increments, it seems hard to find sensible conditions on $X$ under which Lemma~\ref{lem:sup_LP_ladder_height} still holds. Moreover, it is also not clear if a statement similar to Remark~\ref{rem:LP_below_log} is valid. In view of these observations, the restriction that $X$ is a L\'{e}vy process seems quite reasonable.
\end{remark}
\subsection{L\'{e}vy processes}\label{sec:sup_LP_one_sided}
It is not hard to extend Theorem~\ref{thm:LP_at_RW_times} to the case that the inner process is a L\'{e}vy process. We state the result in the next theorem which completes the proof of Theorem~\ref{thm:summary_LP_at_RW_times} announced in the introduction.
\begin{thm}\label{thm:LP_at_Levy_times}
   Let $(X_t)_{t \geq 0}$ and $(Y_t)_{t \geq 0}$ be two independent L\'{e}vy processes such $\E{}{X_1} = 0$, $X_1 \in \mathcal{X}(\alpha)$ and $Y_1 \in \mathcal{X}(\beta)$ for some $\alpha, \beta \in (0,1]$. 
\begin{enumerate}
 \item If $\E{}{Y_1} = 0$, then 
\[
   T^{-1/4} \precsim \pr{}{\sup_{t \in [0,T]} X(\abs{Y_t}) \leq 1} \precsim T^{-1/4} \, (\log T)^{1/4 + 1/(\alpha \wedge \beta)}, \quad T \to \infty.
\]
  \item If $\E{}{Y_1} \neq 0$, then 
\[
   T^{-1/2} \precsim \pr{}{\sup_{t \in [0,T]} X(\abs{Y_t}) \leq 1} \precsim T^{-1/2} \, (\log T)^{1/(\alpha \wedge \beta)}, \quad T \to \infty.
\]
\end{enumerate}
In either case, the lower bound also holds without the assumption of stretched exponential moments.
\end{thm}
\begin{proof}
Upper bound:\\ 
Clearly, we have for all $T >0$ that
\begin{equation}\label{eq:comp_LP_vs_RW}
 \pr{}{\sup_{t \in [0,T]} X(\abs{Y_t}) \leq 1} \leq \pr{}{\sup_{n=1,\dots, \lfloor T \rfloor} X(\abs{Y_n}) \leq 1}.
\end{equation}
Since $Y$ is a L\'{e}vy process, $(Y_n)_{n \geq 1} = \left( \sum_{k=1}^n (Y_k - Y_{k-1}) \right)_{n \geq 1} \stackrel{d}{=} (S_n)_{n \geq 1}$ where $S$ is a random random walk whose increments are equal in distribution to $Y_1$. In particular, the assumptions of Theorem~\ref{thm:LP_at_RW_times} are fulfilled proving the upper bound in both cases.\\
Lower bound for the case $\E{}{Y_1} = 0$: \\
Again, we have that
\begin{align*}
 \pr{}{\sup_{t \in [0,T]} X(\abs{Y_t}) \leq 1} \geq \pr{}{\sup_{t \in [0,c\sqrt{T}]} X_t \leq 1} \, \pr{}{\sup_{t \in [0,T]} \abs{Y_t} \leq c \sqrt{T} }
\end{align*}
Note that by Doob's inequality applied to the submartingale $(\abs{Y_t})_{t \geq 0}$, we obtain that
\[
  \pr{}{\sup_{t \in [0,T]} \abs{Y_t} \leq c\sqrt{T} } = 1 - \pr{}{\sup_{t \in [0,T]} \abs{Y_t} > c\sqrt{T} } \geq 1 - \frac{E{}{Y_T^2}}{c^2 T} = 1 - E{}{Y_1^2} / c^2 = 1/2
\]
for $c := \sqrt{2 \, \E{}{Y_1^2}}$. We have used that $\E{}{Y_t^2} = t \cdot E{}{Y_1^2}$ for a square integrable L\'{e}vy martingale. This proves the lower bound if $\E{}{Y_1} = 0$.\\
Lower bound for the case $\E{}{Y_1} \neq 0$:\\
As before, for any $c > \abs{\E{}{Y_1}}$, we have
\[
 \pr{}{\sup_{t \in [0,T]} X(\abs{Y_t}) \leq 1} \geq \pr{}{\sup_{t \in [0,c T]} X_t \leq 1} \, \pr{}{\sup_{t \in [0,T]} \abs{Y_t} \leq c T }.
\]
Next, since $\abs{Y_t} \leq \abs{Y_t - \E{}{Y_t}} + \abs{\E{}{Y_t}}$ and $\E{}{Y_t} = \E{}{Y_1} \cdot t$ for a L\'{e}vy process, it follows that
\begin{align*}
 \pr{}{\sup_{t \in [0,T]} \abs{Y_t} \leq c T } &\geq \pr{}{\sup_{t \in [0,T]} \abs{Y_t - \E{}{Y_t} } \leq (c - \abs{\E{}{Y_1}} ) T } \\
&\geq 1 - \frac{ \E{}{\abs{Y_T - \E{}{Y_T}}^2} }{ (c- \abs{\E{}{Y_1}})^2 \, T^2  } = 1 - \frac{ \E{}{\abs{Y_1 - \E{}{Y_1}}^2} }{ (c- \abs{\E{}{Y_1}})^2 \, T  } \to 1
\end{align*}
as $T \to \infty$. We have again used Doob's inequality and the fact that $\E{}{\abs{Y_T - \E{}{Y_T}}^2} = \E{}{\abs{Y_1 - \E{}{Y_1}}^2} \cdot T$. This completes the proof of the lower bound.
\end{proof}
\begin{remark}
The above theorem can be strengthend if $X$ is a symmetric L\'{e}vy process and $Y$ is a subordinator. Assume w.l.o.g.\ that $Y_1 \geq 0$ a.s. Then $Z := X \circ Y$ is a symmetric L\'{e}vy process (see e.g.\ Lemma~2.15 of \cite{kyprianou:2006}). In particular,
\[
  \pr{}{\sup_{t \in [0,T]} Z_t \leq 1} \precsim \pr{}{\sup_{n \in [0,\lfloor T \rfloor]} Z_n \leq 1}\asymp T^{-1/2},
\]
without any additional assumption of moments, see e.g.\ Proposition~1.4 of \cite{dembo-gao:2011}. This oberservation suggests that Theorem~\ref{thm:LP_at_RW_times} and \ref{thm:LP_at_Levy_times} remain true under much weaker integrability conditions. In the proof of the upper bound, we needed stretched exponential moments in order to ensure that the distance of $M_{n-1}$ and $M_n$ does not become too large when $M_{n-1} < M_n$. This allowed us (at the cost of a lower order term) to consider the supremum of the process $X$ over the whole interval from $0$ to the maximum of the absolute value of the random walk up to time $N$ instead of the set $\set{\abs{S_1}, \dots, \abs{S_N}}$. Yet, even for a deterministic increasing sequence $(s_n)_{n \geq 1}$ such that $s_N \to \infty$ as $N \to \infty$ and a Brownian motion $(B_t)_{t \geq 0}$, it is not obvious to find conditions on $(s_n)_{n \geq 1}$ such that
\[
   \pr{}{\sup_{n=1,\dots,N} B(s_n) \leq 1} \asymp \pr{}{\sup_{t \in [0,s_N]} B_t \leq 1} \asymp s_N^{-1/2}.
\]
We refer to \cite{aurzada-baumgarten:2011} for related results.
\end{remark}
\section{Two-sided processes}\label{sec:two-sided}
In Sections \ref{sec:sup_self_sim_one_sided}, \ref{sec:sup_RW_one_sided} and \ref{sec:sup_LP_one_sided}, the outer process $X = (X_t)_{t \geq 0}$ had the index set $[0,\infty)$, so it was only possible to evaluate $X$ over the range of the absolute value of the inner process $Y$. In this section, we work with two-sided processes $X = (X_t)_{t \in \mathbb{R}}$ allowing us to consider the one-sided exit problem for the process $X \circ Y$.\\
In Section~\ref{sec:sup_self_sim_two_sided}, we assume that $X$ is a two-sided process with independent branches defined in \eqref{eq:def-two-sided-proc} and that the inner process $Y$ is a self-similar continuous process before turning to the case of random walks and L\'{e}vy processes in Section~\ref{sec:sup_RW_two_sided}. Finally, if $X$ is a fractional Brownian motion indexed by $\mathbb{R}$, the branches of $X$ are not independent (unless $X$ is a two-sided Brownian motion). We provide a brief discussion of this case in Section~\ref{sec:fBm}.
\subsection{Continuous self-similar processes}\label{sec:sup_self_sim_two_sided}
Here we prove a more general version of Theorem~\ref{thm:summary_self_sim_2} which follows from Theorem~\ref{thm:cont_self_sim_proc_two_sided}, Lemma~\ref{lem:int_cond_lower_tail_1} and \ref{lem:inf_sup_expectation}. 
\begin{thm} \label{thm:cont_self_sim_proc_two_sided}
Let $(X_t)_{t \in \mathbb{R}}$ be a two-sided process generated by $X^-$ and $X^+$ with
\[
	\pr{}{\sup_{t \in [0,T]} X_t^- \leq 1} \asymp T^{-\theta^-}, \quad \pr{}{\sup_{t \in [0,T]} X_t^+ \leq 1} \asymp T^{-\theta^+}, \qquad T \to \infty
\]
for some $\theta^-, \theta^+ > 0$. Let $(Y_t)_{t \geq 0}$ denote an independent self-similar process of index $H$ with continuous paths such that $Y_0 = 0$ and
\begin{equation}\label{eq:inf_sup_expectation}
	\E{}{\left(-\inf_{t \in [0,1]} Y_t \right)^{-\theta^-} \, \left(\sup_{t \in [0,1]} Y_t \right)^{-\theta^+} } < \infty.
\end{equation}
Then
\[
 \pr{}{\sup_{t \in [0,T]} X(Y_t) \leq 1} \asymp T^{-H(\theta^- + \theta^+)}, \quad T \to \infty.
\] 
\end{thm}
\begin{proof}
Lower bound: Using the mutual independence of $X^-, X^+$ and $Y$, we get
\begin{align*}
  &\pr{}{\sup_{t \in [0,T]} X(Y_t) \leq 1} \\
&\quad \geq \pr{}{\sup_{t \in [0,T]} \abs{Y_t} \leq T^H} \, \pr{}{\sup_{t \in [0,T^H]} X_t^+ \leq 1} \, \pr{}{\sup_{t \in [0,T^H]} X_t^- \leq 1} \\
&\quad \asymp \pr{}{\sup_{t \in [0,1]} \abs{Y_t} \leq 1} \, T^{-H \theta^-} \, T^{-H \theta^+}.
\end{align*}
In the last step, we have used the self-similarity of $Y$.\\
Upper bound: Denote by $I$ and $M$ the infimum and maximum process of $Y$. By assumption, we can choose a constant $C$ such that for all $T > 0$
\[
	\pr{}{\sup_{t \in [0,T]} X_t^- \leq 1} \leq C \, T^{-\theta^-}, \quad \pr{}{\sup_{t \in [0,T]} X_t^+ \leq 1} \leq C \, T^{-\theta^+}.
\]
 Since the branches $X^+$ and $X^-$ of $X$ are independent, the fact that $Y_0 = 0$ and $Y$ has continuous paths, we have
\begin{align*}
 \pr{}{\sup_{t \in [0,T]} X(Y_t) \leq 1} &= \pr{}{\sup_{t \in [0,-I_T]} X_t^- \leq 1, \sup_{t \in [0,M_T]} X_t^+ \leq 1}\\
&= \E{}{\pr{}{\sup_{t \in [0,-I_T]} X_t^- \leq 1 \vert \mathcal{F}^Y_T } \, \pr{}{\sup_{t \in [0,M_T]} X_t^+ \leq 1 \vert \mathcal{F}^Y_T}  } \\
&\leq C^2 \, \E{}{ (-I_T)^{-\theta^-} \, (M_T)^{-\theta^+} } \\
&=  C^2 \, \E{}{(-I_1)^{-\theta^-} \, M_1^{-\theta^+}} \, T^{-H(\theta^- + \theta^+)}.
\end{align*}
Since the last expectation is finite by assumption, the proof is complete.
\end{proof}
The applicability of Theorem~\ref{thm:cont_self_sim_proc_two_sided} hinges on the verification that the expectation in \eqref{eq:inf_sup_expectation} is finite. The next lemma states such a result. In fact, it turns out that \eqref{eq:inf_sup_expectation} is not harder to verify than the integrability condition of Theorem~\ref{thm:cont_self_sim_proc} if the small deviations of $Y$ satisfy a rather weak condition. 
\begin{lemma}\label{lem:inf_sup_expectation} 
Let $\eta_1, \eta_2 >0$. Assume that 
\[
	\E{}{\left(- \inf_{t \in [0,1]} Y_t\right)^{-\eta_1} } + \E{}{\left(\sup_{t \in [0,1]} Y_t\right)^{-\eta_2} } < \infty.
\]
Moreover, assume that for some $\gamma > 0$, one has
\begin{equation} \label{eq:small_dev_Y}
	\pr{}{\sup_{t \in [0,1]} \abs{Y_t} \leq \epsilon} \precsim \exp\left(-\epsilon^{-\gamma}\right), \qquad \epsilon \downarrow 0. 
\end{equation}
Then the expectation in \eqref{eq:inf_sup_expectation} is finite for any $\theta^- \in (0,\eta_1)$ and $\theta^+ \in (0,\eta_2)$.
\end{lemma}
\begin{proof}
Note that
 \begin{align}
  &\E{}{ ( - I_1)^{-\theta^-} \, M_1^{-\theta^+} } \leq  \E{}{ ( - I_1)^{-\theta^-} \, M_1^{-\theta^+} ; - I_1 \leq \epsilon, M_1 \leq \epsilon } \notag \\
&\qquad + \E{}{ ( - I_1)^{-\theta^-} \, M_1^{-\theta^+} ; - I_1 > \epsilon } + \E{}{ ( - I_1)^{-\theta^-} \, M_1^{-\theta^+} ; M_1 > \epsilon } \notag \\
&\quad \leq \E{}{ ( - I_1)^{-\theta^-} \, M_1^{-\theta^+} ; \sup_{t \in [0,1]} \abs{Y_t} \leq \epsilon } + \epsilon^{-\theta^-} \E{}{M_1^{-\theta^+} } + \epsilon^{-\theta^+} \E{}{(-I_1)^{-\theta^-} } \label{eq:proof_expect_inf_sup_1}.
\end{align}
The two latter expectations are finite due the assumptions on the integrability of $M_1$ and $I_1$. 
Next, for $\epsilon < 1$, we can write
\begin{align*}
\E{}{ ( - I_1)^{-\theta^-} \, M_1^{-\theta^+} ; \sup_{t \in [0,1]} \abs{Y_t} \leq \epsilon } = \sum_{k=1}^\infty \E{}{ (- I_1)^{-\theta^-} \, M_1^{-\theta^+} ; \epsilon^{k+1} < \sup_{t \in [0,1]}  \abs{Y_t} \leq \epsilon^k }. 
\end{align*}
We can choose $p > 1$ such that $p \theta^- < \eta_1$ and $p \theta^+ < \eta_2$. Let $q > 1$ such that $1/p + 1/q = 1$. Using H\"older's inequality in the second estimate, we get
\begin{align}
 &\E{}{ (- I_1)^{-\theta^-} \, M_1^{-\theta^+} ; \epsilon^{k+1} < \sup_{t \in [0,1]}  \abs{Y_t} \leq \epsilon^k } \notag \\
&\quad \leq \epsilon^{-\theta^-(k+1)} \E{}{M_1^{-\theta^+} ; \epsilon^{k+1} < -I_1, \sup_{t \in [0,1]}  \abs{Y_t} \leq \epsilon^k } \notag \\
&\qquad + \epsilon^{-\theta^+(k+1)} \E{}{(-I_1)^{-\theta^-} ; \epsilon^{k+1} < M_1, \sup_{t \in [0,1]}  \abs{Y_t} \leq \epsilon^k } \notag \\
&\quad \leq \epsilon^{-\theta^-(k+1)} \E{}{M_1^{-p\theta^+}}^{1/p} \pr{}{\sup_{t \in [0,1]}  \abs{Y_t} \leq \epsilon^k }^{1/q} \notag \\
&\qquad + \epsilon^{-\theta^+(k+1)} \E{}{(-I_1)^{-p\theta^-}}^{1/p} \pr{}{\sup_{t \in [0,1]}  \abs{Y_t} \leq \epsilon^k }^{1/q}.\label{proof_expect_inf_sup_of_BM3}
\end{align}
By our choice of $p$ and the integrability assumption of the lemma, we see that both expectations in the last expression are finite. 
Next, \eqref{eq:small_dev_Y} implies that there is a constant $C > 0$ such that $\pr{}{\sup_{t \in [0,1]}  \abs{Y_t} \leq \epsilon} \leq C \exp(-\epsilon^{-\gamma})$ for all $\epsilon < 1$. Therefore, for $\epsilon < 1$ and $k \geq 1$, we have for any $\eta > 0$ that
\begin{align*}
  \sum_{k=1}^\infty \epsilon^{-\eta k} \pr{}{\sup_{t \in [0,1]}  \abs{Y_t} \leq \epsilon^k}^{1/q} \leq C  \sum_{k=1}^\infty \epsilon^{-\eta k} \exp\left(- \epsilon^{-\gamma k}/q \right) < \infty.
\end{align*}
Hence, in view of \eqref{proof_expect_inf_sup_of_BM3}, it follows that all expressions in \eqref{eq:proof_expect_inf_sup_1} are finite.
\end{proof}
\begin{remark}\label{rem:inf_sup_expec_BM}
In view of Lemma~\ref{lem:int_cond_lower_tail_1}, one can easily check whether the assumptions of Lemma~\ref{lem:inf_sup_expectation} are fulfilled. For instance, if $Y=B$ is a Brownian motion, then 
\[
\pr{}{-\inf_{t \in [0,1]} B_t \leq \epsilon} = \pr{}{\sup_{t \in [0,1]} B_t \leq \epsilon} = \pr{}{\abs{B_1} \leq \epsilon} \sim \sqrt{2/\pi} \, \epsilon, \quad \epsilon \downarrow 0   
\]
Hence, 
\[
   \E{}{\left(-\inf_{t \in [0,1]} B_t\right)^{-\theta^-} \left(\sup_{t \in [0,1]} B_t \right)^{-\theta^+} } < \infty, \quad \theta^-,\theta^+ \in (0,1).
\]
Note that we cannot use H\"{o}lder's inequality to establish this result if $\theta^-,\theta^+ \geq 1/2$.
\end{remark}
We can now state a result for iterated Brownian motion (cf.\ \cite{burdzy:1993}).
\begin{cor}
Let $(B_t)_{t \in \mathbb{R}}$ be a two-sided Brownian motion and $(W_t)_{t \geq 0}$ denote another independent Brownian motion. Then
\[
 \pr{}{\sup_{t \in [0,T]} B(W_t) \leq 1} \asymp T^{-1/2}, \quad T \to \infty.
\] 
\end{cor}
\begin{proof}
This follows directly from Theorem~\ref{thm:cont_self_sim_proc_two_sided} and Remark~\ref{rem:inf_sup_expec_BM}. 
\end{proof}
Of course, we can apply Theorem~\ref{thm:cont_self_sim_proc_two_sided} to any two-sided process $X$ whose branches have survival exponents strictly smaller than one and $Y$ being an independent Brownian motion. Examples for $X$ therefore include two-sided intergrated Brownian motion (survival exponent $\theta^+ = \theta^- = 1/4$, two-sided symmetric L\'{e}vy processes ($\theta^+ = \theta^- = 1/2$) and fBm ($\theta^+ = \theta^- = 1-H$ where $H$ is the Hurst parameter of the fBm (here, one has to use an obvious extension of Theorem~\ref{thm:cont_self_sim_proc_two_sided} taking into account that the survial probability decays like $T^{-(1-H) + o(1)}$)). Of coure, the branches of $X^-$ and $X^+$  need not have the same distribution. 
\subsection{Two-sided L\'{e}vy processes at random walk or L\'{e}vy times}\label{sec:sup_RW_two_sided}
Let us now consider the one-sided exit problem for the process $(X(S_n))_{n \geq 0}$ where $S$ is again a random walk and $X$ is a two-sided L\'{e}vy process, i.e.\ the branches of $X$ are independent L\'{e}vy processes. The next theorem shows that the survival exponent is $1/2$ under suitable integrability conditions regardless of the sign of $\E{}{S_1}$ in contrast to Theorem~\ref{thm:LP_at_RW_times}. 
\begin{thm}\label{thm:two_sided_LP_at_RW_times}
   Let $(X_t)_{t \in \mathbb{R}}$ denote a two-sided L\'{e}vy process with branches $X^+, X^-$, $\E{}{X_1^-} = \E{}{X_1^+} = 0$ and $X_1^-, X_1^+ \in \mathcal{X}(\alpha)$ for some $\alpha \in (0,1]$. Let $(Y_n)_{n \geq 1}$ denote a sequence of i.i.d.\ random variables independent of $X$ with $Y_1 \in \mathcal{X}(\beta)$ for some $\beta \in (0,1]$.  Let $S_n = Y_1 + \dots + Y_n$. Then
\[
 \pr{}{\sup_{n=1,\dots,N} X(S_n) \leq 1} = N^{-1/2 + o(1)}, \quad N \to \infty.
\]
More specifically: 
\begin{enumerate}
 \item If $\E{}{Y_1} = 0$, then 
  \[
    N^{-1/2} \precsim \pr{}{\sup_{n=1,\dots,N} X(S_n) \leq 1} \precsim N^{-1/2} \, \, (\log N)^{1/2 + 2/(\alpha \wedge \beta)}, \quad N \to \infty.
  \]
  \item 
If $\E{}{Y_1} \neq 0$, then
  \[
    N^{-1/2} \precsim \pr{}{\sup_{n=1,\dots,N} X(S_n) \leq 1} \precsim N^{-1/2} \, (\log N)^{1/(\alpha \wedge \beta)}, \quad N \to \infty.
  \]
\end{enumerate}
In either case, the lower bound also holds without the assumption of stretched exponential moments on $X_{-1}$, $X_1$ and $Y_1$.
\end{thm}
\begin{proof}
 The lower bound can be established as in the proof of Theorem~\ref{thm:LP_at_Levy_times} if $\E{}{Y_1} = 0$. If $\E{}{Y_1} > 0$ (say), using that $\inf_{n \geq 1} S_n$ is a finite random variable a.s., the result follows along similar lines.
\\
For the upper bound, assume first that $\E{}{Y_1} > 0$. Then
\begin{align*}
 p_N &:= \pr{}{\sup_{n=1,\dots,N} X(S_n) \leq 1} \leq \pr{}{\sup_{n=1,\dots,N} X^+(M_n) \leq 1} \\
&\leq \pr{}{M_N \leq c N} + \pr{}{M_N \geq c N, \sup_{n=1,\dots,N} X^+(M_n) \leq 1}.
\end{align*}
The second line is just equation \eqref{eq:upper_bound_RW_pos_exp}, so the result follows from the proof the upper bound of Theorem~\ref{thm:LP_at_RW_times}.\\
If $\E{}{Y_1} < 0$, the result follows by applying the theorem to the random walk $(-S_n)_{n \geq 1}$. \\
Let us finally consider the case $\E{}{Y_1} = 0$. Let $f(N) := \sqrt{\log N}$, $N \geq 1$. Note that
\begin{align*}
 p_N &\leq \pr{}{\sup_{n=1,\dots,N} \abs{S_n} \leq N^{1/2}/f(N)  } + \\
&\quad + \pr{}{ M_N \leq N^{1/2}/f(N), -I_N > N^{1/2}/f(N), \sup_{n=1,\dots,N} X(S_n) \leq 1 } \\
&\quad + \pr{}{ M_N > N^{1/2}/f(N), -I_N \leq N^{1/2}/f(N), \sup_{n=1,\dots,N} X(S_n) \leq 1 } \\
&\quad + \pr{}{ M_N > N^{1/2}/f(N), -I_N > N^{1/2}/f(N), \sup_{n=1,\dots,N} X(S_n) \leq 1 } \\
&=: J_1(N) + J_2(N) + J_3(N) + J_4(N).
\end{align*}
First, recall that $J_1(N) = o(N^{-1/2})$ (cf. \eqref{eq:de-acosta_result}). It remains to estimate the terms $J_2$ and $J_4$. The term $J_3$ can be dealt with analogously to $J_2$.\\
\emph{Step 1:}\\
\begin{align*}
 J_2(N) &\leq \pr{}{M_N \leq N^{1/4}, -I_N > N^{1/2}/f(N), \sup_{n=1,\dots,N} X^-(-I_n) \leq 1} \\
&\quad + \pr{}{N^{1/4} \leq M_N \leq N^{1/2}/f(N), -I_N > N^{1/2}/f(N), \sup_{n=1,\dots,N} X(S_n) \leq 1} \\
&=: K_{2,1}(N) + K_{2,2}(N).
\end{align*} 
Let us now find upper bounds for $K_{2,j}$ for $j=1,2$. Denote by $\sigma^+(n)$ resp.\ $\sigma^-(n)$ the $n$-th time that the random walk $S$ reaches a new maximum resp.\ minimum and by $\mathcal{H}^+_n := S_{\sigma^+(n)} - S_{\sigma^+(n-1)}$ resp.\ $\mathcal{H}^-_n := -(S_{\sigma^-(n)} - S_{\sigma^-(n-1)})$ the corresponding ascending resp.\ descending ladder heights. \\
\emph{Step 2:}\\
First, note that
\begin{align*}
 K_{2,1}(N) = \E{}{\indic{M_N \leq N^{1/4}} \,\indic{ -I_N >N^{1/2}/f(N)} \, \pr{}{ \sup_{n=1,\dots,N} X^-(-I_n) \leq 1 \vert \mathcal{F}_N^S } }.
\end{align*}
Next, proceeding just as in the proof of Theorem~\ref{thm:LP_at_RW_times}, we obtain for $N$ large enough that
\begin{align}
 q_N &:= \indic{ -I_N > N^{1/2}/f(N)} \, \pr{}{ \sup_{n=1,\dots,N} X^-(-I_n) \leq 1 \vert \mathcal{F}_N^S } \notag \\
&\leq \indic{ -I_N > N^{1/2}/f(N)} \pr{}{ \sup_{t \in [0,-I_N]} X^-_t \leq 1 + c (\log N)^\rho \vert \mathcal{F}_N^S } \notag \\
&\quad + \pr{}{ \bigcup_{n : \sigma^-(n) \leq N} \set{ \sup_{t \in [-S_{\sigma^{-}(n-1)},-S_{\sigma^{-}(n)}]} X^-_t - X^-(S_{\sigma^-(n-1)}) \geq c(\log N)^\rho }  \vert \mathcal{F}_N^S } \notag \\
&\leq C \frac{c(\log N)^\rho}{N^{1/4}/\sqrt{f(N)}} + \sum_{n=1}^N \pr{}{\sup_{t \in [0,\mathcal{H}_n^-]} X^-_t > c(\log N)^\rho \vert \mathcal{F}_N^S}. \label{eq:proof_q_N}
\end{align}
In the last inequality, we have used Remark~\ref{rem:LP_below_log} assuming that $\rho \geq 1/\alpha$. Using that the $\mathcal{H}^-_n$ are i.i.d., this shows that
\begin{align*}
 K_{2,1}(N) \leq  C \frac{c(\log N)^{\rho + 1/4}}{N^{1/4}} \, \pr{}{M_N \leq N^{1/4}} + N \pr{}{\sup_{t \in [0,\mathcal{H}_1^-]} X^-_t > c(\log N)^\rho}.
\end{align*}
Applying the second part of Lemma~\ref{lem:one-sided-exit-normalized-RW} with $\tilde{f}(N) := N^{1/4}$ to the first summand and Lemma~\ref{lem:sup_LP_ladder_height} to the second, we obtain with $\rho :=1/(\alpha \wedge \beta)$ for $c$ large enough
\begin{equation}\label{eq:K_2_1}
 K_{2,1}(N) \precsim \frac{c(\log N)^{1/(\alpha \wedge \beta) + 1/4}}{N^{1/4}} \, N^{-1/4} + N^{-1/2} \precsim (\log N)^{1/(\alpha \wedge \beta) + 1/4} \, N^{-1/2}.
\end{equation}
Let us now find an upper bound on $K_{2,2}$. Set
\[
 r_N := \indic{N^{1/4} \leq M_N \leq N^{1/2}/f(N)}  \pr{}{ \sup_{n=1,\dots,N} X^+(M_n) \leq 1 \vert \mathcal{F}_N^S }.
\]
Since $X^-$ and $X^+$ are independent, we have in view of \eqref{eq:proof_q_N} and $r_N \leq 1$ that
\begin{align}
 K_{2,2} \leq \E{}{r_N q_N } \leq \frac{d(\log N)^{1/(\alpha \wedge \beta) + 1/4}}{N^{1/4}} \, \E{}{r_N } + N \pr{}{\sup_{t \in [0,\mathcal{H}_1^-]} X^-_t > c(\log N)^{1/(\alpha \wedge \beta)}} \label{eq:K_2_2}.
\end{align}
Let $a(k) := \sum_{l=1}^k 2^{-(l+1)} = (1-2^{-k})/2$, $k \geq 1$. Since $a(N) \to 1/2$, we can find $\gamma(N)$ such that $N^{a(\gamma(N))} \geq N^{1/2}/f(N)$. Indeed, this just amounts to
\begin{align}
   a(\gamma(N)) = \frac{(1 - 2^{-\gamma(N)})}{2} \geq \frac{\log (N^{1/2}/f(N))}{\log N}= \frac{1}{2} - \frac{\log \log N}{2 \log N},
\end{align}
i.e.\
\[
   \gamma(N) \geq \frac{1}{\log 2} \, \log \left( \frac{\log N}{\log \log N} \right).
\]
Hence, it suffices to set $\gamma(N) := \lceil (\log \log N)/\log 2 \rceil$. \\
Next, note that
$\set{N^{1/4} \leq M_N \leq N^{1/2}/f(N)} \subseteq \set{N^{a(1)} \leq M_N \leq N^{a(\gamma(N))} }$ and proceeding as in \eqref{eq:sup_discr_to_cont}, we obtain 
\begin{align*}
r_N &\leq \sum_{k=1}^{\gamma(N)-1} \indic{N^{a(k)} \leq M_N \leq N^{a(k+1)}}  \pr{}{ \sup_{n=1,\dots,N} X^+(M_n) \leq 1 \vert \mathcal{F}_N^S } \\
&\leq \sum_{k=1}^{\gamma(N)-1} \indic{N^{a(k)} \leq M_N \leq N^{a(k+1)}}  \pr{}{ \sup_{t \in [0,M_N]} X^+_t \leq 1 + c (\log N)^{1/(\alpha \wedge \beta)} \vert \mathcal{F}_N^S } \\
&\quad + \sum_{k=1}^{\gamma(N)-1} \indic{N^{a(k)} \leq M_N \leq N^{a(k+1)}} \sum_{n=1}^N \pr{}{\sup_{t \in [0,\mathcal{H}_n^+]} X^+_t > c (\log N)^{1/(\alpha \wedge \beta)} \vert \mathcal{F}_N^S} \\
&\leq \sum_{k=1}^{\gamma(N)-1} \indic{M_N \leq N^{a(k+1)}}  \pr{}{ \sup_{t \in [0,N^{a(k)}]} X^+_t \leq 2c (\log N)^{1/(\alpha \wedge \beta)} } \\
&\quad + \gamma(N) \sum_{n=1}^N \pr{}{\sup_{t \in [0,\mathcal{H}_n^+]} X_t > c (\log N)^{1/(\alpha \wedge \beta)} \vert \mathcal{F}_N^S}.
\end{align*}
Taking expectations and keeping in mind Lemma~\ref{lem:sup_LP_ladder_height} and Remark~\ref{rem:LP_below_log}, we conclude that
\begin{align*}
\E{}{r_N} &\leq \sum_{k=1}^{\gamma(N)-1} \pr{}{ M_N \leq N^{a(k+1)} }  \pr{}{ \sup_{t \in [0,N^{a(k)}]} X^+_t \leq 2c (\log N)^{1/(\alpha \wedge \beta)} } + o(N^{-1/2})\\
&\leq C \, \sum_{k=1}^{\gamma(N) - 1} \pr{}{ M_N \leq N^{a(k+1)} } (\log N)^{1/(\alpha \wedge \beta)} \, N^{-a(k) / 2}  + o(N^{-1/2})
\end{align*}
for some $c$ large enough. In view of Lemma~\ref{lem:one-sided-exit-normalized-RW} (second part), we can find constants $C_1$ and $N_0$ such that for $N \geq N_0$
\[
 \pr{}{ M_N \leq N^{a(k+1)} } \leq C_1 \, N^{a(k+1) -1/2 }, \qquad k = 1,2,\dots
\]
(Note that we can get such a uniform estimate, see Remark~\ref{rem:sawyer}). Hence, for $N$ large enough, we obtain
\begin{align*}
 &\sum_{k=1}^{\gamma(N)-1} \pr{}{ M_N \leq N^{a(k+1)} } (\log N)^{1/(\alpha \wedge \beta)} \, N^{-a(k) / 2} \\
&\quad \leq C_1 (\log N)^{1/(\alpha \wedge \beta)} \, \sum_{k=1}^{\gamma(N)-1} N^{a(k+1) - 1/2 - a(k) / 2} = (\gamma(N) - 1) \, C_1 \, (\log N)^{1/(\alpha \wedge \beta)} \, N^{-1/4}
\end{align*}
since $a(k+1) - a(k) / 2 = 1/4$. This shows that $\E{}{r_N} \precsim \gamma(N) \, (\log N)^{1/(\alpha \wedge \beta)} \, N^{-1/4}$ and therefore, we deduce from \eqref{eq:K_2_2} and the definition of $\gamma(N)$ that
\[
  K_{2,2}(N) \precsim (\log \log N) \, (\log N)^{2/(\alpha \wedge \beta) + 1/4} \, N^{-1/2}, \quad N \to \infty.
\]
Combining this with \eqref{eq:K_2_1}, it follows that
\begin{equation}\label{eq:J_2_two_sided}
   J_2(N) \precsim (\log \log N) \, (\log N)^{2/(\alpha \wedge \beta) + 1/4} \, N^{-1/2}, \quad N \to \infty.
\end{equation}
\emph{Step 3}: \\
It remains to consider $J_4$. The line of reasoning is now clear. Setting $g(N) := N^{1/2}/f(N) = \sqrt{ N/\log N}$, one can check that $J_4$ is bounded from above by
\begin{align*}
 &\E{}{\indic{M_N \geq g(N)} \pr{}{ \sup_{n=1,\dots,N} X^+_{M_n} \leq 1 \vert \mathcal{F}_N^S } \, \indic{-I_N \geq g(N)} \pr{}{\sup_{n=1,\dots,N} X^-_{-I_n} \leq 1 \vert \mathcal{F}_N^S } } \\
&\precsim \pr{}{\sup_{t \in [0,g(N)]} X^+_t \leq c (\log N)^{1/(\alpha \wedge \beta)} } \, \pr{}{\sup_{t \in [0,g(N)]} X^-_t \leq c (\log N)^{1/(\alpha \wedge \beta)} } + N^{-1/2} \\
&\precsim \left((\log N)^{1/(\alpha \wedge \beta)} g(N)^{-1/2}\right)^2 = (\log N)^{1/2 + 2/(\alpha \wedge \beta)} \, N^{-1/2}.
\end{align*}
This finishes the proof.
 \end{proof}
 \begin{remark}\label{rem:two-sided-LP}
 The proof reveals that the survival exponent is equal to $1/2$ no matter if $\E{}{Y_1} = 0$ or not for quite different reasons. If $\E{}{Y_1} > 0$, $S_N / N \to \E{}{Y_1}$ by the law of large numbers, so the random walk diverges to $+\infty$ with speed $N$ and the survival probability is determined by the right branch $X^+$ of $X$. \\
 If $\E{}{Y_1} = 0$, the random walks oscillates and typical fluctuations are of order $\pm \sqrt{N}$. The survival probability up to time $N$ is therefore approximately equal to the probability that both $X^+$ and $X^-$ stay below $1$ until time $\sqrt{N}$. By independence of $X^+$ and $X^-$, this probability is equal to the product of these two probabilities which are each of order $N^{-1/4}$.
 \end{remark}
Clearly, the analogue of Theorem~\ref{thm:LP_at_Levy_times} also holds for two-sided L\'{e}vy processes. We state this result without proof.
\begin{thm}
  Let $(X_t)_{t \in \mathbb{R}}$ denote a two-sided L\'{e}vy process with branches $X^+, X^-$, $\E{}{X_1^-} = \E{}{X_1^+} = 0$ and $X_1^-, X_1^+ \in \mathcal{X}(\alpha)$ for some $\alpha \in (0,1]$. Let $(Y_t)_{t \geq 1}$ be another L\'{e}vy process independent of $X$ with $Y_1 \in \mathcal{X}(\beta)$ for some $\beta \in (0,1]$. 
\begin{enumerate}
 \item If $\E{}{Y_1} = 0$, then for any $\epsilon > 0$, we have
\[
   T^{-1/2} \precsim \pr{}{\sup_{t \in [0,T]} X(Y_t) \leq 1} \precsim T^{-1/2} \, (\log T)^{1/2 + 1/(\alpha \wedge \beta)}, \quad T \to \infty.
\]
  \item If $\E{}{Y_1} \neq 0$, then 
\[
   T^{-1/2} \precsim \pr{}{\sup_{t \in [0,T]} X(Y_t) \leq 1} \precsim T^{-1/2} \, (\log T)^{1/(\alpha \wedge \beta)}, \quad T \to \infty.
\]
\end{enumerate}
\end{thm}

\subsection{Fractional Brownian motion}\label{sec:fBm}
Let $(X_t)_{t \in \mathbb{R}}$ denote a fBm with Hurst parameter $H \in (0,1)$, i.e.\ $X$ is a centered Gaussian process with covariance  
\[
 \E{}{X_t X_s} = \frac{1}{2} \, \left( \abs{t}^{2H} + \abs{s}^{2H} - \abs{t-s}^{2H} \right), \qquad s,t \in \mathbb{R}.
\]
If $s< 0 <t$, one can check that $\E{}{X_t X_s} > 0$ if $H < 1/2$ and $\E{}{X_t X_s} < 0$ if $H > 1/2$. Hence, the branches of a fBm are not independent unless $H=1/2$ and Theorem~\ref{thm:cont_self_sim_proc_two_sided} is not applicable. However, it is not difficult to find an appropriate generalization if the survival exponent of the two-sided process is known. We now state such a result for fBm.
\begin{prop}
 Let $(X_t)_{t \in \mathbb{R}}$ denote a fBm with Hurst parameter $H \in (0,1)$ and $(Y_t)_{t \geq 0}$ a self-similar process of index $\lambda >0$ with continuous paths. Assume that for any $0 < \eta < 1$, it holds that
\[
 \E{}{\left(\sup_{t \in [0,1]} Y_t \right)^{-\eta}} < \infty, \quad \E{}{\left(-\inf_{t \in [0,1]} Y_t \right)^{-\eta}} < \infty.
\]
Then
\[
 \pr{}{\sup_{t \in [0,T]} X(Y_t) \leq 1} = T^{-\lambda + o(1)}, \quad T \to \infty.
\] 
In particular, the survival exponent does not depend on $H$.
\end{prop}
\begin{proof}
By Theorem~3 of \cite{molchan:1999a}, we have for any $H \in (0,1)$ that
\begin{align*}
 \pr{}{\sup_{t \in [-T,T]} X_t \leq 1} = T^{-1 + o(1)}, \quad T \to \infty.
\end{align*}
The lower bound of the proposition can be proved just as in Theorem~\ref{thm:cont_self_sim_proc_two_sided}. \\
For the upper bound, fix $\epsilon \in (0,1)$. Then we can find a constant $C > 0$ such that $\pr{}{\sup_{t \in [-T,T]} X_t \leq 1} \leq C T^{-1 + \epsilon}$ for any $T > 0$. Moreover, since the paths of $Y$ are continuous, we have that
\begin{align*}
 &\pr{}{\sup_{t \in [0,T]} X(Y_t) \leq  1} = \E{}{ \pr{}{\sup_{t \in [I_T, M_T]} X_t \leq 1 \vert \mathcal{F}^Y_T} } \\
& \leq \E{}{ \indic{M_T \geq -I_T} \pr{}{\sup_{[I_T, -I_T]} X_t \leq 1 \vert \mathcal{F}^Y_T} } +  \E{}{ \indic{M_T \leq -I_T} \pr{}{\sup_{[-M_T, M_T]} X_t \leq 1 \vert \mathcal{F}^Y_T} } \\
&\leq C \left( \E{}{(-I_T)^{-1 + \epsilon}} + \E{}{M_T^{-1+ \epsilon}} \right) = C \left( \E{}{(-I_1)^{-1 + \epsilon}} + \E{}{M_1^{-1+ \epsilon}} \right) \, T^{-\lambda + \lambda \epsilon}. 
\end{align*}
\end{proof}

\bibliographystyle{abbrvnat}
\bibliography{global_bib}
\end{document}